\documentclass{amsart}

\usepackage{graphicx,epsfig}
\usepackage[dvips]{psfrag}
\usepackage{enumerate}
\usepackage{amssymb}

\newtheorem{theo}{Theorem}
\newtheorem{Coro}{Corollary}
\newtheorem{Ques}{Question}

\newtheorem{theor}{Theorem}[section]
\newtheorem{lemm}[theor]{Lemma}

\newtheorem{coro}[theor]{Corollary}

\newtheorem{defi}[theor]{Definition}
\newtheorem{prop}[theor]{Proposition}

\newtheorem{rema}[theor]{Remark}

\newtheorem{ques}[theor]{Question}

\addtocounter{figure}{0}

\numberwithin{equation}{section}
\numberwithin{figure}{section}

\newcommand{\eqdef}{\stackrel{\scriptscriptstyle\rm def}{=}}

 \def\BB{{\mathbb B}} \def\CC{{\mathbb C}}
  
\def\EE{{\mathbb E}} \def\FF{{\mathbb F}} \def\GG{{\mathbb G}}

  \def\KK{{\mathbb K}}

\def\LL{{\mathbb L}}

 \def\RR{{\mathbb R}} \def\SS{{\mathbb S}}
\def\TT{{\mathbb T}}

\def\YY{{\mathbb Y}}

\def\si{\sigma}
\def\Si{\Sigma}

\def\De{\Delta}

\def\Ga{\Gamma}

\def\ga{\gamma}

\def\ve{\varepsilon}

    \def\cT{{\mathcal
T}}
\def\cC{{\mathcal C}}    \def\cU{{\mathcal
U}}

\def\cY{{\mathcal Y}}

\def\loc{{\operatorname{loc}}}

\def\Diff{{\operatorname{Diff}}^1}

\def\Difr{{\operatorname{Diff}}^r}

\begin{document}

\title{Fragile cycles}

\author{Ch. Bonatti}
\address{Institut de Math\'ematiques de Bourgogne, BP 47 870, 1078 Dijon Cedex, France}
\email{bonatti@u-bourgogne.fr}

\author{L. J. D\'\i az}
\address{Depto. Matem\'atica PUC-Rio, Marqu\^es de S. Vicente 225
22453-900 Rio de Janeiro RJ  Brazil} \email{lodiaz@mat.puc-rio.br}


\date{\today}


\begin{abstract}
We study diffeomorphisms $f$ with heterodimensional cycles, that is, 
 heteroclinic cycles associated to saddles $p$ and $q$ with different indices. Such a cycle is called fragile if there is no diffeomorphism close to $f$ with a robust cycle associated to hyperbolic sets containing the continuations of $p$ and $q$.
We construct a codimension one submanifold of $\Difr(\SS^2\times \SS^1)$, $r\ge 1$, that consists of diffeomorphisms with fragile heterodimensional cycles.  Our construction holds for any manifold of dimension $\ge 4$. 
\end{abstract}
\maketitle



\section{Introduction}\label{s.introduction}
In the late sixties, Newhouse constructed the first examples of
$C^2$-open sets of non-hyperbolic surface diffeomorphisms. Any such set
 $\cU$ consists of diffeomorphisms  with \emph{$C^2$-robust
homoclinic tangencies}:
every diffeomorphism   $f\in \cU$ has a hyperbolic set
$K_f$ (depending continuously on $f$) whose stable and unstable
manifolds have non-transverse intersections, see \cite{N68}.

Later,  in \cite{N79}, Newhouse proved that homoclinic tangencies of
surface diffeomorphisms  can be {\it{stabilized\/}}: 
given a diffeomorphism $f$ with
a homoclinic tangency associated to a saddle $p_f$, 
there is a $C^2$-open set whose closure contains $f$ 
and which consists of diffeomorphisms $g$ with robust homoclinic tangencies
associated to hyperbolic sets
$K_g$ containing the continuation $p_g$ of $p_f$. 
In particular, these results show that homoclinic tangencies always generate
$C^2$-robust homoclinic tangencies. In fact, robust homoclinic
tangencies are present in all known examples of open sets of
non-hyperbolic surface diffeomorphisms. Let us observe that
 homoclinic tangencies of $C^1$-diffeomorphisms
defined on surfaces cannot be stabilized, see \cite{Gugu}.

Similarly, all known examples of $C^1$-open sets formed by
non-hyperbolic diffeomorphisms 
exhibit $C^1$-robust {\emph{heterodimensional cycles,}} that is,
cycles relating the invariant manifolds of two hyperbolic sets of
different $s$-indices (dimension of the stable bundle). Note that
the existence of such cycles can only occur in dimension $\ge 3$.

We wonder if, as in the case of homoclinic tangencies of
$C^2$-diffeomorphisms, heterodimensional cycles can  be made
$C^1$-robust and can be $C^1$-stabilized. A first partial answer to this
question is given in \cite{BD08}: heterodimensional cycles associated to
periodic saddles whose indices differ by one 
generate (by arbitrarily small $C^1$-perturbations) 
$C^1$-robust heterodimensional cycles.  
In some extend, the results in \cite{BD08} are a version of the ones
by Newhouse in \cite{N68,N79} in the context of  $C^1$-heterodimensional
cycles. 

However, compared with Newhouse's results for homoclinic tangencies, the ones in
\cite{BD08} have an important disadvantage: While the hyperbolic
sets with the robust homoclinic tangencies in \cite{N68,N79}
contain continuation of the saddle with the initial tangency,
the hyperbolic sets involved in robust cycles in \cite{BD08} do in general not contain the continuations of the saddles in the initial cycle. However this precise question can be important  
for understanding the global dynamics of non-hyperbolic
diffeomorphisms.
Let us discuss more this question in more detail.

Following Conley theory \cite{Con}
and motivated by spectral decomposition theorems \cite{Newhouse,A03},
 this global dynamics  is structured using
{\emph{homoclinic}} or/and \emph{chain recurrence classes} as
``elementary" pieces of dynamics, see the 
definitions below. 
One aims 
to describe the dynamics of each piece and the relations between
different pieces (cycles), for further details see \cite[Chapter
10.3-4]{BDVbook} and \cite{bible}. 

In general, 
the homoclinic class of a hyperbolic periodic
point is contained in its chain recurrence class.
 An important property is that  for  $C^1$-generic diffeomorphisms 
 homoclinic classes
and chain recurrence classes  of periodic points coincide,
\cite[Remarque 1.10]{BC04}. However,
in non-generic situations, two different homoclinic classes (even of saddles of different indices) may be joined by  a cycle,
hence they are contained in the same chain recurrence class. 
A question is when  one can join them in a $C^1$-robust way by small
perturbations. This occurs if the cycle can be
stabilized. For instance, this is specially important for understanding the indices of the periodic points in an elementary piece of dynamics.

While the above explains why the stabilization of cycles is relevant,
let us now provide the precise 
definitions of the concepts involved. First, recall that a hyperbolic basic set $K$ of a
diffeomorphism $f$ has a (uniquely defined) {\emph{continuation}}
$K_g$ for all $g$ close to $f$: $K_g$ is a hyperbolic basic set, close to $K$, and the dynamics of $f|_{K}$ and $g|_{K_g}$ are conjugate. The {\emph{$s$-index}} of a hyperbolic transitive set is the dimension of its stable bundle.

\begin{defi}[Robust continuations of cycles]
{\em{$\,$

{$\bullet$} The diffeomorphism $f$ is said to have a {\emph{heterodimensional
cycle}} associated to hyperbolic basic sets
 $K$ and $L$ if these sets have different $s$-indices and their
stable and unstable manifolds meet cyclically:
$$
W^s(K,f)\cap  W^u(L,f)\ne\emptyset \quad \mbox{and} \quad
W^u(K,f)\cap W^s(L,f)\ne\emptyset.
$$

{$\bullet$} The  cycle associated to $K$ and $L$ is {\emph{
$C^1$-robust}} if there is a $C^1$-neighborhood $\cU$ of $f$ such
that for all $g\in \cU$ the hyperbolic continuations $K_g$ and
$L_g$ of $K$ and $L$ have a heterodimensional cycle.

{$\bullet$} A heterodimensional cycle associated to a pair of
saddles $p$ and $q$ of $f$ {\emph{can be $C^1$-stabilized}} if
every $C^1$-neighborhood $\cU$ of $f$ contains a diffeomorphism
$g$ with hyperbolic basic sets $K_g\ni p_g$ and $L_g\ni q_g$ having a robust heterodimensional cycle. Here $p_g$ and $q_g$ are the continuations of
$p$ and $q$ for $g$.}}
\end{defi}

\begin{defi}[Fragile cycle]{\em{
A heterodimensional cycle associated to a pair of saddles is
{\emph{$C^1$-fragile}} if it cannot be $C^1$-stabilized.}}
\end{defi}

The previous discussion leads to the following 
question that we address in this paper:
{\emph{
Can every heterodimensional cycle be
$C^1$-stabilized?}}

\smallskip

The results in \cite{BDK}
give a positive answer to this question for ``most" types of
{\emph{coindex one}}
 heterodimensional cycles, that is,  related to saddles whose $s$-indices differ by one. Indeed in \cite{BDK} it is  proved that fragile  coindex one cycles  associated to saddles $p$ and $q$ exhibit a quite specific geometry:
\begin{itemize}
\item
The homoclinic classes of $p$ and $q$ are both trivial.
\item
The  {\emph{central  eigenvalues}} of $p$ and $q$ are all real and
positive.\footnote{The definition of central eigenvalues is a little
intricate. Assuming that the $s$-index of $p$ is bigger than the
one of $q$, the central eigenvalues correspond to the weakest
contracting direction of $p$ and the weakest expanding direction
of $q$.}
\item
There is a well-defined one-dimensional  orientable central bundle
$E^c$ along the cycle (i.e. defined on some closed set containing
the saddles $p$ and $q$ in the cycle and a pair of heteroclinic
orbits $x\in W^s(p)\cap W^u(q)$ and
 $y\in W^u(p)\cap W^s(q)$), but the cycle diffeomorphism does not preserve
the orientation of $E^c$. The cycle is \emph{twisted} by the terminology in 
\cite{ASY06}.
\end{itemize}

In this paper we provide examples of fragile coindex one cycles, see
Theorem~\ref{t.fragil}. 
 
\subsection{Definitions and statement of results}
Recall that the {\emph{homoclinic class}} of a hyperbolic periodic
point $p$, denoted by $H(p,f)$, is the closure of the transverse
intersections of the stable and unstable manifolds of the orbit of
$p$. The homoclinic class $H(p,f)$ coincides with the closure of the set of all saddles
$q$ {\emph{homoclinically related with $p$,}} i.e. the stable
manifold of the orbit of $q$ transversely meets the unstable
manifold of the orbit of $p$ and vice-versa. A homoclinic class is
{\emph{non-trivial}} if it contains at least two different orbits.

Let us now recall the definition of a {\emph{chain recurrence class.}} 
A finite sequence of points $(x_i)_{i=0}^n$ is an
{\emph{$\epsilon$-pseudo-orbit}} of a diffeomorphism $f$ if
$\mbox{dist}(f(x_i),x_{i+1})<\epsilon$ for all $i=0,\dots,n-1$. A
point $x$ is {\emph{chain recurrent}} for $f$ if
 for every $\epsilon>0$ there is an $\epsilon$-pseudo-orbit
 $(x_i)_{i=0}^n$, $n\ge 1$,
starting and ending at $x$ (i.e. $x=x_0=x_n$). The {\emph{chain recurrent set}} $R(f)$ of $f$ is the set of all chain recurrent points. This
set splits into disjoint {\emph{chain recurrence classes}}: the
class $C(x,f)$ of $x\in R(f)$ is the set of points $y$ such that
for every $\epsilon>0$ there are $\epsilon$-pseudo-orbits joining
$x$ to $y$ and $y$ to $x$. 
A periodic point $p$ of  $f$ is {\emph{isolated}} if its {\emph{chain
recurrence class}} coincides with its orbits. In this case, the
orbit of $p$ is the maximal invariant set in some \emph{filtrating
neighborhood}. This implies that the homoclinic class of $p$
is $C^1$-robustly trivial (i.e. the homoclinic class of  $p_g$ is
trivial for every $g$ close to $f$).

We are now ready to state our main result.

\begin{theo}\label{t.fragil}
There is an open set $\cU$ of $\Diff (\SS^2\times \SS^1)$ and a
codimension one submanifold $\Sigma$ contained in $\cU$ with the
following property: For every $f\in \cU$ there are  hyperbolic
saddles $p_f$ and $q_f$ with
different $s$-indices depending continuously on $f$  such that
\begin{enumerate}
\item every  $f\in \Sigma$ has
a heterodimensional cycle associated to $p_f$ and $q_f$,
\item the set $\cU\setminus\Sigma$ is the union of  two connected sets $\, \cU^+$ and
$\,\cU^-$ such that
\begin{itemize}
\item
for every $f\in \cU^+$ the saddle $p_f$ is isolated,
\item
for every $f\in \cU^-$ the saddle $q_f$ is isolated.
\end{itemize}
\end{enumerate}
\end{theo}

Note that if two hyperbolic basic sets $K_f$ and $L_f$ have a
heterodimensional cycle then the chain recurrence classes of any
pair of  saddles $a_f\in K_f$ and $b_f\in L_f$ coincide. In
particular, if the cycle associated to $K_f$ and $L_f$ is robust
then the chain recurrence classes of $a_g$ and $b_g$ are the same
for all $g$ in some neighborhood of $f$. In particular, the chain recurrence
class $C(a_g,g)=C(b_g,g)$ is non-trivial and the saddles $a_g$ and
$b_g$ are both non-isolated for $g$. Thus
Theorem~\ref{t.fragil} implies that the cycles in $\Si$ cannot be made
robust. This implies  the following:

\begin{Coro}
\label{c.fragil} The submanifold $\Si$ in Theorem~\ref{t.fragil}
consists of diffeomorphisms having  $C^1$-fragile
heterodimensional cycles.
\end{Coro}

As for any $n>3$ the set $\SS^2\times \SS^1$ can be embedded as a normally contracting manifold  in a ball $\BB^n$, we obtain the
following.

\begin{Coro}
Any compact manifold $M$ with $\dim M>3$ supports diffeomorphisms
with $C^1$-fragile cycles.
\end{Coro}

As our examples demand a somewhat specific topological
configuration, the following question arises naturally.

\begin{Ques} Does every $3$-manifold admit diffeomorphisms with $C^1$-fragile
heterodimensional cycles?
\end{Ques}

The examples presented in this paper display many interesting and somehow unexpected properties. There are also many important aspects of their dynamics yet unexplored. Thus, after completing our construction, in Section~\ref{s.conclusion} we conclude with a discussion about the properties of our 
examples.

\medskip

This paper is organized as follows. 
In the first step of our construction, in Section~\ref{s.auxiliary}, 
we build an auxiliary  Morse-Smale
vector field $X$ on the $3$-sphere $\SS^3$.
In Section~\ref{s.surgery}, we consider  a surgery 
in $\SS^3$
(associated to some identifications by
a local diffeomorphism $\Psi$ of $\SS^3$). This surgery provides
 a diffeomorphism $F_\Psi$ defined on $\SS^2 \times \SS^1$ 
induced by the time-one map $F_0=X_1$ of the vector field $X$
and the gluing map $\Psi$. We also see how the dynamics of $F_\Psi$
depends on the gluing map $\Psi$. In Section~\ref{ss.neighborhood}, we study the dynamics of
diffeomorphisms close to $F_\Psi$. Finally, in Section~\ref{s.psi}, we
choose the gluing map $\Psi$ to get  a diffeomorphism $F_\Psi$ with a fragile cycle 
and construct the submanifold $\Si$ consisting of diffeomorphisms with
fragile cycles. The paper is closed with a discussion section.



\section{An auxiliary vector field on
$\SS^3$}\label{s.auxiliary}

In this section we construct a Morse-Smale vector field
defined on the three sphere $\SS^3$ 
whose non-wandering set consists of 
singular points.
This vector field also
satisfies some normally hyperbolic properties. We now go to the
details of this construction.

We consider the sphere $\SS^3$  as  the union of
two solid tori $\cT_1$ and $\cT_2$ with the same boundary
$\partial \cT_1=\partial \cT_2=\TT^2$. A simple closed curve of
$\TT^2$ is a \emph{$\cT_i$-meridian} if it is not $0$-homotopic in
$\TT^2$ but is $0$-homotopic in $\cT_i$.
We consider an  identification of the boundaries of these solid tori
that does not preserve the meridians: the
$\cT_2$-meridians are isotopic in $\cT_1$ to the ``central circle"
of $\cT_1$ and are classically called \emph{$\cT_1$-parallels}.
Similarly, $\cT_1$-meridians are $\cT_2$-parallels.

\subsection{An auxiliary Morse-Smale vector field  $X$ in  $\SS^3$}\label{ss.auxiliary}

Consider a Morse-Smale vector field $X$ defined on $\SS^3$ such
that (see Figure~\ref{f.MorseSmale}):

\begin{figure}[htb]
       \psfrag{s1}{$s_1$}
\psfrag{s2}{$s_2$}
\psfrag{p1}{$p_1$}
\psfrag{p2}{$p_2$}
\psfrag{r1}{$r_1$}
\psfrag{r2}{$r_2$}
\psfrag{q1}{$q_1$}
\psfrag{q2}{$q_2$}
\psfrag{su}{$\sigma^u$}
\psfrag{ss}{$\sigma^s$}
\psfrag{x}{$X$}
\psfrag{g1s}{$\gamma_1^s$}
\psfrag{g2s}{$\gamma_2^s$}
\psfrag{T2}{$\cT_2$}
\psfrag{T1}{$\cT_1$}
\psfrag{T}{$\TT$}
\includegraphics[width=7.5cm]{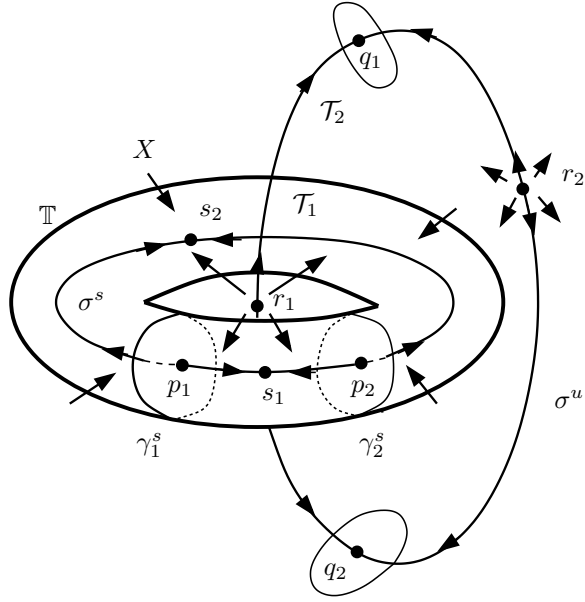}\hspace{1cm}
    \label{f.MorseSmale}

 \caption{The Morse-Smale vector field $X$ defined on $\SS^3$.}

\end{figure}

\begin{enumerate}
\item
$X$ is transverse to $\partial \cT_1=\partial \cT_2=\TT^2$.
\item
The solid torus $\cT_1$ is attracting and the solid torus $\cT_2$
is repelling: The positive orbit of any point $x\in \TT^2$ enters
(and remains) in the interior of $\cT_1$ and its negative orbits
enters (and remains) in the interior of $\cT_2$.
\item
\label{i.circles}
The maximal invariant set of $X$ in $\cT_1$ is a normally
hyperbolic (contracting) circle $\sigma^s$. Analogously, the
maximal invariant set of $X$ in $\cT_2$ is a normally hyperbolic
(repelling) circle $\sigma^u$.
\item The limit set of $X$ is
$\{s_1,s_2,r_1,r_2,p_1,p_2,q_1,q_2\}$, where
\begin{itemize}
\item
$s_1,s_2\in \sigma^s$ are attracting singularities,
\item
$r_1,r_2\in \sigma^u$ are repelling singularities,
\item
$p_1,p_2\in \sigma^s$ are saddle singularities of $s$-index $2$,
and
\item
$q_1,q_2\in \sigma^u$ are saddle singularities of $s$-index $1$.
\end{itemize}
\item
\label{i.Xsepartices} The two (one-dimensional) separatrices
of the unstable manifold of the singularity $p_i$,
$i=1,2$, are  contained in
$W^s(s_1)$ and in   $W^u(s_2)$. A similar
assertion holds for the separatrices of the stable manifold of
$q_i$, $i=1,2$, that are contained in $W^u(r_1)$ and $W^u(r_2)$.
\item The local stable manifold of  $p_i$, $i=1,2$, is a 
$2$-disk
contained in $\cT_1$  whose boundary 
\[
\partial W^s_{loc}(p_i) =
W^s_{loc}(p_i)\cap
\partial\cT_1\eqdef \ga_i^s
\]
is a  $\cT_1$-meridian. Similarly, 
the local unstable manifold of $q_i$, $i=1,2$,  is a  $2$-disk
contained in $\cT_2$ whose boundary
\[
\partial
W^u_{loc}(q_i)=W^u_{loc}(q_i)\cap
\partial \cT_2\eqdef \ga_i^u,
\]
 is a  $\cT_2$-meridian.
\item
\label{i.points} For every $i,j\in\{1,2\}$, the curve $\gamma_i^s$
is transverse to $\gamma_j^u$ and the intersection $\ga_i^s\cap
\ga_j^u$ is  exactly one point $x_{ij}$.
\end{enumerate}

\begin{figure}[htb]

\psfrag{s1}{$s_1$}
\psfrag{s2}{$s_2$}
\psfrag{p1}{$p_1$}
\psfrag{p2}{$p_2$}
\psfrag{ss}{$\sigma^s$}
\psfrag{ss+}{$\sigma^s_+$}
\psfrag{ss-}{$\sigma^s_-$}
\psfrag{su+}{$\sigma^u_+$}
\psfrag{su-}{$\sigma^u_-$}
\psfrag{g1s}{$\gamma^s_1$}
\psfrag{g2s}{$\gamma^s_2$}
\psfrag{g1u}{$\gamma^u_1$}
\psfrag{gu2}{$\gamma^u_2$}
\psfrag{b1s}{$b^s_1$}
\psfrag{b2s}{$b^s_2$}
\psfrag{C1s}{$C_1^s$}
\psfrag{C2s}{$C_2^s$}
\psfrag{b1u}{$b^u_1$}
\psfrag{b2u}{$b^u_2$}
\psfrag{R}{$R$}
\psfrag{T1}{$\cT_1$}
\psfrag{a1s}{$a^s_1$}
\psfrag{a2s}{$a^s_2$}

   \includegraphics[width=11cm]{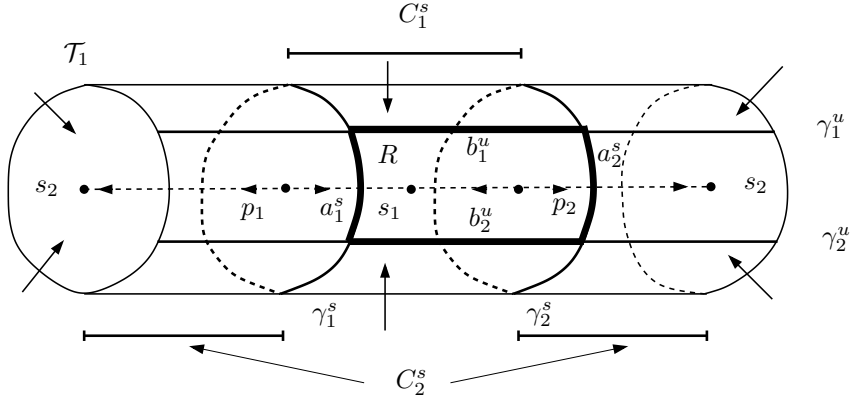}\hspace{1cm}
    \caption{The rectangle $R$.}
    \label{f.rectangle}
\end{figure}

\begin{rema}
[Dynamics of the vector field $X$] \label{r.stablecylinders} $\,$
\em{
\begin{enumerate}

\item \label{i.boundary}
The boundary  $\TT^2$ of the solid torus $\cT_1$ is the union of
two cylinders $\CC_1^s$ and $\CC_2^s$ with disjoint interiors and
the same boundary $\ga_1^s \cup \ga_2^s$. The notation is chosen
such that $W^s(s_i)\cap \TT^2$ is the interior of the cylinder
$\CC_i^s$.

Similarly,   $\TT^2=\partial \cT_2$ is the union of the cylinders
$\CC_1^u$ and $\CC_2^u$  bounded by $\ga_1^u$ and $\gamma_2^u$ and
whose interiors are the intersections $\TT^2\cap W^u(r_1)$ and
$\TT^2\cap W^u(r_2)$, respectively. See Figure~\ref{f.rectangle}.
\item \label{i.cilinders}
As a consequence of item (\ref{i.points}) in the definition of the vector field $X$, the intersection
$\overline{\CC_1^u} \cap \overline{\CC_1^s}$ is a ``rectangle"
$R$ such that
$$
R= \overline{W^s(s_1)\cap W^u(r_1)\cap
\TT^2}
$$ 
and its boundary $\partial R$ of $R$ is the union of four curves
$a_1^s,\, a_2^s,\,b_1^u, \, b_2^u$ with disjoint interiors such
that
\begin{equation}\label{e.boundaryR}
a_1^s\subset \ga_1^s,  \quad a_2^s\subset \gamma_2,\quad
b_1^u\subset \gamma^u_1,\quad
 b_2^u\subset \gamma^u_2.
\end{equation}
See Figure~\ref{f.rectangle}.
Note that the interiors  of $b_1^u$ and $b_2^u$ are contained in
$W^s(s_1)$ and the interiors of  $a_1^s$ and $a_2^s$ are contained
in $W^u(r_1)$. 
\end{enumerate}
}
\end{rema}

\subsection{Partially hyperbolicity of $X$}\label{ss.partially}

We also assume that the vector field $X$ satisfies the following partially hyperbolic
conditions:

\begin{enumerate}
\item[(1)] There is a partially hyperbolic splitting of  $X$ over the
circle $\sigma^s$ (recall \eqref{i.circles} in Section~\ref{ss.auxiliary})
of the form
$$
T_{\sigma^s}\SS^3=    E^{ss}\oplus E^{cs},
$$
where $E^{cs}$ is a $2$-dimensional central bundle containing the
$X$ direction, and $E^{ss}$ is a strong stable bundle that is
oriented along the circle $\sigma^s$.

A similar condition holds for the circle $\sigma^u$: There is a
partially hyperbolic splitting of  $X$ over $\sigma^u$
of the form
$$
T_{\sigma^u}\SS^3= E^{cu} \oplus E^{uu},
$$
where $E^{cu}$ is a $2$-dimensional center bundle
containing the $X$-direction
 and $E^{uu}$ is
a strong unstable bundle that is oriented along $\sigma^u$.

\item[(2)]
Consider  the two-dimensional strong stable manifold $W^{ss}(\sigma^s)$ of $\sigma^s$ that is tangent to
$X\oplus E^{ss}$ along $\sigma^s$. Define the
local strong stable
manifold of $\sigma^s$ by   $W^{ss}_\loc (\sigma^s)=W^{ss} (\sigma^s) \cap \cT_1$. 
Then the intersection between $W^{ss}_\loc (\sigma^s)$
and  $\TT^2$ consists of two disjoint
$\cT_1$-parallels $\si^s_-$ and $\si^s_+$. Similarly, the intersection
$W^{uu}_\loc (\sigma^u)\cap
\TT^2$ is the disjoint union of  two $\cT_2$-parallels $\si^u_+$ and
$\si^u_-$.
We require that
$$
\CC_1^u \cap \big( \si^s_- \cup \si^s_+ \big)=\emptyset  \quad
\mbox{and} \quad \CC_1^s \cap \big( \si^u_- \cup \si^u_+
\big)=\emptyset.
$$
\end{enumerate}

\begin{figure}[htb]

\psfrag{s1}{$s_1$}
\psfrag{s2}{$s_2$}
\psfrag{p1}{$p_1$}
\psfrag{p2}{$p_2$}
\psfrag{q1}{$q_1$}
\psfrag{q2}{$q_2$}
\psfrag{ss}{$\sigma^s$}
\psfrag{ss+}{$\sigma^s_+$}
\psfrag{ss-}{$\sigma^s_-$}
\psfrag{su+}{$\sigma^u_+$}
\psfrag{su-}{$\sigma^u_-$}
\psfrag{g1s}{$\gamma^s_1$}
\psfrag{g2s}{$\gamma^s_2$}
\psfrag{g1u}{$\gamma^u_1$}
\psfrag{gu2}{$\gamma^u_2$}
\psfrag{b1s}{$b^s_1$}
\psfrag{b2s}{$b^s_2$}
\psfrag{b1u}{$b^u_1$}
\psfrag{b2u}{$b^u_2$}
\psfrag{R}{$R$}
\psfrag{T1}{$\cT_1$}
\psfrag{a1s}{$a^s_1$}
\psfrag{a2s}{$a^s_2$}

   \includegraphics[width=12cm]{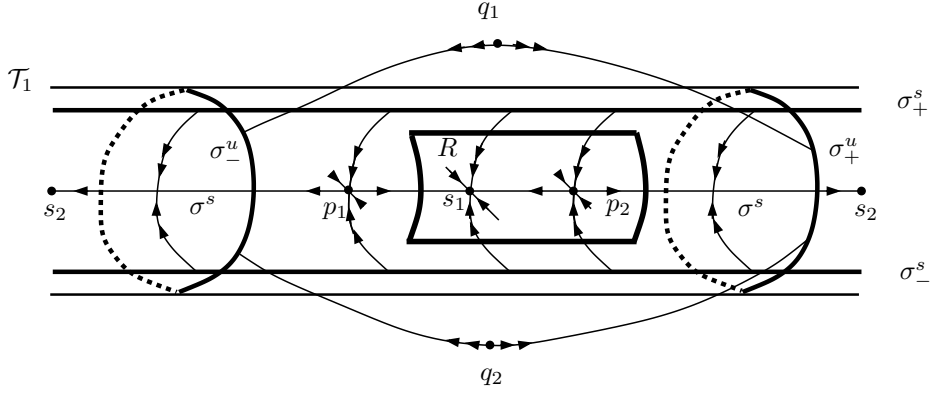}\hspace{1cm}
    \caption{The rectangle $R$ and the ``strong'' manifolds.}
    \label{f.strong}
\end{figure}

Let us explain how this property can be obtained. Recall that
$W^u(r_1)\cap \TT^2$ is the interior of the cylinder $\CC_1^u$
bounded by $\gamma_1^u$ and
$\gamma_2^u$. Thus, since $\gamma_1^u$, $\gamma_2^u$,
$\sigma_-^s,$ and $\sigma_+^s$ are $\cT_1$-parallels (or
equivalently $\cT_2$-meridians), we can assume that $\CC_1^u \cap
\big( \si^s_- \cup \si^s_+ \big)=\emptyset$. See
Figure~\ref{f.strong}. The condition for the cylinder $\CC_1^s$
and the circles $\si^u_+$ and $\si^u_-$ follows identically
noting that $\gamma_1^s$, $\gamma_2^s$, $\sigma_-^u$, and
$\sigma_+^u$ are $\cT_2$-parallels.

By the partially hyperbolic conditions,  the strong stable
manifolds $W^{ss}(p_i)$, $i=1,2$, (tangent to $E^{ss}$ at $p_i$) are well defined and
has dimension one. Similarly, the strong unstable manifold
$W^{uu}(q_1)$ and $W^{uu}(q_2)$ are well  defined and have dimension
one. As a consequence of item (2) above (see Figure~\ref{f.strong}) we have
the following:

\begin{equation}\label{e.boundaryRbis}
\begin{split}
& a_1^s\cap W^{ss}(p_1)=\emptyset, \quad a_2^s\cap
W^{ss}(p_2)=\emptyset,\\ & b_1^u\cap W^{uu}(q_1) =\emptyset, \quad
b_2^u\cap W^u(q_2) =\emptyset.
\end{split}
\end{equation}

\subsection{Transverse heteroclinic intersection}
\begin{figure}[htb]
\psfrag{s1}{$s_1$}
\psfrag{p1}{$p_1$}
\psfrag{p2}{$p_2$}
\psfrag{q1}{$q_1$}
\psfrag{q2}{$q_2$}
\psfrag{G11}{$\Gamma_{1,1}$}
\psfrag{G12}{$\Gamma_{1,2}$}
\psfrag{G21}{$\Gamma_{2,1}$}
\psfrag{G22}{$\Gamma_{2,2}$}

   \includegraphics[width=6cm]{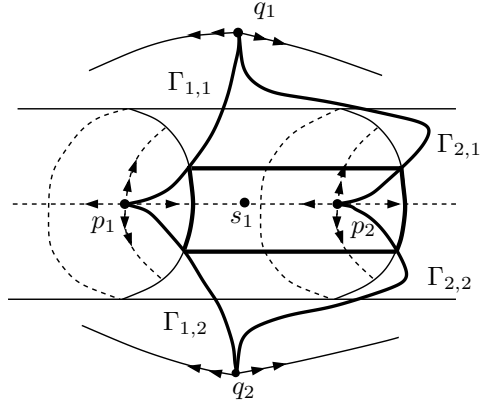}\hspace{1cm}
    \caption{The curves $\Ga_{i,j}$. Global Dynamics}
    \label{f.gaij}
\end{figure}

\begin{figure}[htb]
\psfrag{s1}{$s_1$}
\psfrag{p1}{$p_1$}
\psfrag{p2}{$p_2$}
\psfrag{q1}{$q_1$}
\psfrag{q2}{$q_2$}
\psfrag{G11}{$\Gamma_{1,1}$}
\psfrag{G12}{$\Gamma_{1,2}$}
\psfrag{G21}{$\Gamma_{2,1}$}
\psfrag{G22}{$\Gamma_{2,2}$}

   \includegraphics[width=6cm]{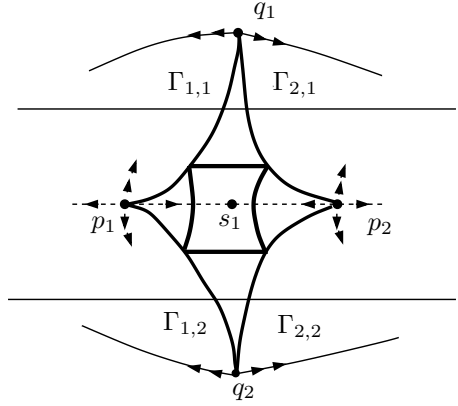}\hspace{1cm}
    \caption{Outline of the curves $\Ga_{i,j}$.}
    \label{f.gaijbis}
\end{figure}

Consider the ``corner'' points of the rectangle $R$,
\begin{equation}\label{e.xij}
a_1^s\cap b_1^u\eqdef x_{1,1}, \quad a_1^s \cap b_2^u \eqdef
x_{1,2}, \quad a_2^s\cap b_1^u\eqdef x_{2,1}, \quad a_2^s \cap
b_2^u\eqdef x_{2,2}.
\end{equation}
By definition, the $\omega$ and $\alpha$-limits of the point $x_{i,j}$ are
the singularities $p_i$ and  $q_j$, respectively.
Denote by $\Ga_{i,j}$ the closure of the orbit of $x_{i,j}$ (see
Figures~\ref{f.gaij} and \ref{f.gaijbis}). As the intersection between $W^s(p_i)$ and
$W^u(q_j)$ is exactly the orbit of $x_{i,j}$, we have
\begin{equation}\label{e.curvasGa}
\Ga_{i,j}\eqdef \overline{W^s(p_i)\cap W^u(q_j)}=\{p_i\} \cup
\{q_j\} \cup \big(W^s(p_i) \pitchfork W^u(q_j) \big) .
\end{equation}

\begin{rema}\label{r.Gaij}
{\em{The curve $\Ga_{i,j}$ is a $C^1$-invariant normally
hyperbolic compact segment.}} 
\end{rema}

This remark is a standard consequence of the following facts:

\begin{itemize}
\item
The point $x_{i,j}$ is a transverse heteroclinic intersection associated
to the singularities $p_i$ and $q_j$.
\item
The partial hyperbolicity hypothesis at the singularities
$p_i$ and $q_j$
 implies that  $W^{ss}(p_i)$ and  
$W^{uu}(q_j)$ are well defined. 
\item
By construction, recall equation~\eqref{e.boundaryRbis}, 
 $x_{i,j}\not\in W^{ss}(p_i) \cup
W^{uu}(q_j)$.
\item
The curve $\Ga_{i,j}$ is the closure of the orbit of $x_{i,j}$.
\end{itemize}

\subsection{Invariant manifolds of the segments $\Ga_{i,j}$}\label{ss.invariant}

\begin{figure}[htb]

\psfrag{p1}{$p_1$}
\psfrag{p2}{$p_2$}
\psfrag{q1}{$q_1$}
\psfrag{G11}{$\Ga_{1,1}$}
\psfrag{g21}{$\Ga_{2,1}$}
\psfrag{Wup1}{$W^u(p_1)$}
\psfrag{Wup2}{$W^u(p_2)$}
\psfrag{wuuq1}{$W^{uu}(q_1)$}
\psfrag{wu+q1}{$W^{u,+}(q_1)$}
\psfrag{wu-q1}{$W^{u,-}(q_1)$}
\psfrag{s1}{$s_1$}
\psfrag{b1u}{$b^{u}_1$}
\psfrag{x11}{$x_{1,1}$}
\psfrag{x21}{$x_{2,1}$}

   \includegraphics[width=7cm]{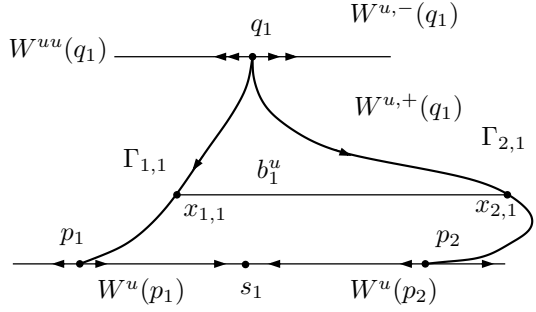}\hspace{1cm}
    \caption{The unstable manifold of $\Ga_{i,j}$.}    \label{f.unstable}
\end{figure}

For each singularity $q_j$, we have that $\big( W^u(q_j)\setminus
W^{uu}(q_j) \big)$ is the disjoint union of two connected
invariant surfaces
$$
\big( W^u(q_j)\setminus  W^{uu}(q_j) \big) \eqdef W^{u,+}(q_j)\cup
W^{u,-} (q_j),
$$
where  $W^{u,+}(q_j)$ contains the interior of the
curves $\Ga_{1,j}$ and $\Ga_{2,j}$. With this notation, the
invariant manifolds of the curve $\Ga_{i,j}$ are
\[
\begin{split}
&W^{u}(\Ga_{i,j})= W^u(p_i) \cup W^{u,+}(q_j) \cup
W^{uu}(q_j),
\\
&W^{s}(\Ga_{i,j})= W^s(q_j) \cup W^{s,+}(p_i) \cup
W^{ss}(p_i).
\end{split}
\]
See Figures~\ref{f.gaijbis} and \ref{f.unstable}. Note that these manifolds are
injective $C^1$-immersions of $[0,1]\times \RR$.

Figures~\ref{f.gaij}, \ref{f.gaijbis}, and \ref{f.unstable} suggest that the curves $\Ga_{1,j}$
and $\Ga_{2,j}$ form a ``cusp" at the point $q_j$. This geometric configuration 
will play a key role in our construction. So let us define precisely
what we mean by a {\emph{cusp.}}

A topological two-disk $\GG$ contained in the interior of a smooth surface $S$
{\emph{has a cusp}} at a point $p\in\partial \GG$ if for every
$\varepsilon
>0$ there are a convex cone $C_\varepsilon$ of angle $\varepsilon$
at $p$ and a neighborhood $U_\varepsilon$ of $p$ in $S$ such that
$\GG\cap U_\varepsilon\subset C_\varepsilon$. Given a curve
$\gamma$ contained in the interior of the surface $S$, a point $q$ in the interior of 
$\gamma$ is a {\emph{cusp}} of
$\gamma$ if there is a topological disk $\FF$ whose boundary
contains $\gamma$ and has  a cusp at  $q$.

\begin{rema} \label{r.curvesGaij}{\em{
Recall that the interior of the curves  $\Ga_{1,i}$ and
$\Ga_{2,i}$ are disjoint from $W^{uu}(q_i)$. Moreover, the
interior of these curves  are the orbits of the points $x_{1,i}$
and $x_{2,i}$, respectively. These two curves are connected by the
segment $b_i^u\subset W^u(q_i)$ that is disjoint from $W^{uu}(q_i)$,
recall \eqref{e.boundaryR} and see Figure~\ref{f.unstable}.  The partial hyperbolicity hypothesis
now implies that $\Ga_{1,i}$ and $\Ga_{2,i}$ are ``central curves''
arriving to $q_1$ from the same side of $W^{uu}(q_i)$. These two
conditions imply that 
\begin{equation}
\label{e.Gaqi}
\Ga (q_i)\eqdef \Ga_{1,i}\cap \Ga_{2,i}, \quad i=1,2,
\end{equation}
is a curve with a \emph{cusp singularity} at $q_1$, see Figures~\ref{f.gaij},
\ref{f.gaijbis}, 
 and
\ref{f.unstable}.}}
\end{rema}

The unstable manifold $W^u(\Ga(q_i))$ of $\Ga(q_i)$ is the set
$W^u(\Ga_{1,i})\cup W^u(\Ga_{2,i})$. Noting that the interior of
the ``strip" $W^u(\Ga(q_i))$ is $W^{u,+}(q_i)$, $i=1,2$,  we get
\begin{equation}
\label{e.q} W^u(\Ga (q_i))=  W^u(p_1)\cup W^u(p_2) \cup
W^{u,+}(q_i) \cup W^{uu}(q_i).
\end{equation}
We observe that the set $W^u(\Ga (q_i))$ is an injective
$C^1$-immersion of a connected surface with boundary. Equivalent
statements hold for
\begin{equation}
\label{e.Gapi}
\Ga(p_i)\eqdef \Ga_{i,1}\cup \Ga_{i,2}, \quad i=1,2,
\end{equation}
and its stable manifold
\begin{equation}
\label{e.p} W^s(\Ga(p_i))= W^s(q_1)\cup W^s(q_2) \cup W^{s,+}(p_i)
\cup W^{ss}(p_i),
\end{equation}
where $W^{s,+}(p_i)$ is the  component of
 $\big( W^s(p_i)\setminus
W^{ss}(p_i) \big)$ containing the interior of the curves
$\Ga_{i,1}$ and $\Ga_{i,2}$.

With this notation, the sides $a_i^s$ and $b_j^u$
of the rectangle $R$ (see \eqref{e.boundaryR}) satisfy the
following property
\begin{equation} \label{e.sides}
 a_i^s \subset W^s(\Ga(p_i))\cap \TT^2 \quad
\mbox{and}\quad b_j^u\subset W^u(\Ga(q_j))\cap \TT^2, \qquad
i,j\in\{1,2\}.
\end{equation}


\subsection{Central bundles}
\label{ss.central}

\begin{figure}[htb]

\psfrag{p1}{$p_i$}
\psfrag{q2}{$q_2$}
\psfrag{q1}{$q_1$}
\psfrag{wuG11}{$W^u(\Ga_{i,1})$}
\psfrag{wuG12}{$W^u(\Ga_{i,2})$}
\psfrag{wssp1}{$W^{ss}(p_i)$}
\psfrag{G11}{$\Ga_{i,1}$}
\psfrag{G12}{$\Ga_{i,2}$}
\psfrag{Ec+}{$E^{cs}_+$}
\psfrag{Ec}{$E^{cs}$}

 \includegraphics[width=7cm]{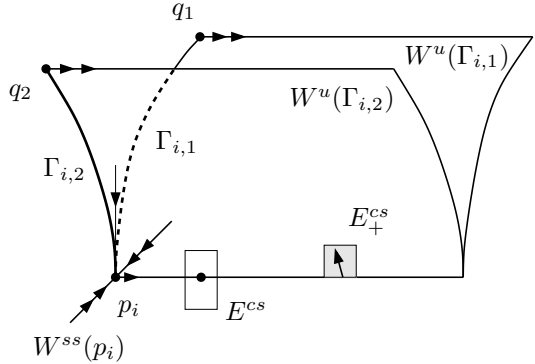}\hspace{1cm}
   
\caption{The open ``half-planes" $E^{cu}_{+}$.}
 \label{f.ecplus}
\end{figure}

In this section, we see that the unstable manifolds of $\Ga_{i,1}$ and 
$\Ga_{i,2}$ touch each other 
at $W^u(p_i)$ tangentially ``coming from the same side" of 
$W^u(p_i)$, see Figure~\ref{f.ecplus}. In the following we will precise 
what this means.

\begin{lemm}[Center stable/unstable
bundles]\label{l.centerbundles}
Given any singularity $p$ of
saddle type with a strong stable direction $W^{ss}(p)$ (tangent
to a strong stable bundle $E^{ss}(p)$) there is a unique invariant
``central" bundle $E^c$ defined over the unstable manifold $W^u(p)$ of $p$ that is transverse at
$p$ to the bundle $E^{ss}$ and has codimension $\dim (E^{ss})$.

A similar property holds for saddle singularities $q$  with a
strong unstable manifold tangent to some strong unstable bundle
$E^{uu}$. In this case, there is a central bundle
defined over
$W^s(q)$ that is
transverse to $E^{uu}$ at $q$ and has codimension $\dim (E^{uu})$.
\end{lemm}

\begin{proof}
We can assume that for every point $x\in W^u_{\loc}(p)$ 
there 
 is defined a
negatively invariant cone-field $\cC^{ss}$ around the strong
stable direction $E^{ss}$. Consider the complement $\cC^c$ of
$\cC^{ss}$. Given $y\in W^u(p)$ there is $t(y)>0$ such that
$X_{-t}(y) \in W^u_{\loc}(p)$ for all $t\ge t(y)$.
Given $y\in  W^u(p)$ it is enough to define
$$
E^{c}(y)\eqdef \{ v\, \colon \, D_y X_{-t}(v) \in
\cC^{c}(X_{-t}(y)) \quad \mbox{for all $t\ge t(y)$}\}.
$$
By construction, the  bundle $E^c(y)$ is transverse to $E^{ss}$ and its dimension
is the codimension of $E^{ss}$. 
This completes the proof of the lemma.
\end{proof}

\begin{rema}\label{r.notation}
{\em{
Applying Lemma~\ref{l.centerbundles} to the singularity $p_i$, we get a two
dimensional bundle $E^c(p_i)$, $i=1,2$, coinciding with the
bundle $E^{cs}$ defined along  the curve $\sigma^s$ in
Section~\ref{ss.partially}. Analogously, the bundle $E^c$ defined
along $W^s(q_j)$ coincides with the bundle $E^{cu}$ along
$\sigma^u$.}}
\end{rema}

\begin{lemm}\label{l.tangent}
The surfaces $W^u(\Ga(q_1))$ and $W^u(\Ga(q_2))$ are tangent to the bundle
$E^{cs}$ along
their intersection $W^u(p_1)\cup W^u(p_2)$. Similarly, the
surfaces $W^s(\Ga(p_1))$ and $W^s(\Ga(p_2))$ are tangent 
to the bundle $E^{cu}$ 
along their intersection $W^s(q_1)\cup
W^s(q_2)$.
\end{lemm}

\begin{proof}
The boundary part of $W^u(\Ga(q_j))$ has three components,
$W^u(p_1)$, $W^u(p_2)$, and $W^{uu}(q_j)$,   
recall  equation \eqref{e.q} and see Figure~\ref{f.unstable}. Moreover, the
surface $W^u(\Ga(q_j))$ is transverse to $W^{ss}(p_i)$. The
uniqueness of the central bundle $E^c$ in
Lemma~\ref{l.centerbundles} implies that for each 
$x\in W^u(p_i)\subset \sigma^s$
the fiber $E^c(x)=E^{cs}(x)$ (recall Remark~\ref{r.notation}) is the tangent space $T_x(W^u(\Ga(q_j)))$. This
implies the lemma.
\end{proof}

\begin{rema}[The open ``half-planes" $E^{cs}_{+}$ and $E^{cu}_{+}$]
\label{r.semiplanes} {\em{ The normally hyperbolic curves
$\Ga_{i,1}$ and $\Ga_{i,2}$ are ``central curves" contained in
$W^s(p_i)$ arriving  to $p_i$ from the same side of $W^s(p_i)\setminus
W^{ss}(p_i)$.
Furthermore, the boundary surfaces  $W^u(\Ga_{i,1})$ and
$W^u(\Ga_{i,2})$ are tangent to $E^{cs}$ along $W^u(p_i)$.  This
implies that, for every $x\in W^u(p_i)$, the vectors in $E^{cs}(x)$
pointing to the interior of $W^u(\Ga(q_1))$) form an open half-plane
$E^{cs}_+(x)$. This half-plane coincides with the vectors of $E^{cs}(x)$
pointing  to the interior of 
 $W^u(\Ga_{i,1})$ or (equivalently)
of
$W^u(\Ga_{i,2})$. See Figure~\ref{f.ecplus}.

For points $y\in W^s(q_j)$, we similarly define the half-plane
$E^{cu}_+(y)$ as the vectors in $E^{cu}$ pointing to the interior of
$W^s(\Ga(p_i))$, $i=1,2$.}}
\end{rema}

\subsection{Position of the invariant manifolds in the basins of $r_1$ and
$s_1$}\label{s.positionofinvariant} Consider a ``small" two-sphere
$\SS^s$ contained in the interior of the solid torus $\cT_1$ that
is transverse to the vector field $X$ and bounds a three-ball
$\BB^s\subset W^s(s_1) \cap \cT_1$ whose interior contains the
singularity $s_1$. Let
\begin{equation}\label{e.etas}
\eta^s\eqdef \SS^s \cap W^{ss}(\si^s).
\end{equation}
Note that $\eta^s$ is a circle that contains the points
\begin{equation}\label{e.yu12}
y_1^u\eqdef W^u(p_1) \cap \SS^s \quad \mbox{and} \quad y_2^u\eqdef
W^u(p_2) \cap \SS^s, \qquad y_1^u,y^u_2\in \eta^s.
\end{equation}
We similarly define a ``small" two-sphere $\SS^u \subset \cT_2$
transverse to $X$ bounding a three-ball $\BB^u\subset W^u(r_1)$ whose
interior contains $r_1$. We define the
circle
$$
\eta^u\eqdef
\SS^u\cap 
 W^{uu}(\si^u) 
$$
and the intersection
points
\begin{equation}\label{e.ys12}
y_1^{s}\eqdef  W^s(q_1) \cap \SS^u \quad \mbox{and} \quad
y^s_2= W^s(q_2) \cap \SS^u, \qquad y_1^s,y^s_2\in \eta^u.
\end{equation}

\begin{rema}\label{r.time}
{\em{Choosing the balls $\BB^s$ and $\BB^u$ small enough, we can assume that
the minimum time that a point takes to go from $\BB^u$ to $\BB^s$ is arbitrarily large. In particular, this time is bigger than $10$:
$X_t (\BB^u) \cap \BB^s=\emptyset$ for all $t\in [0,10]$. }}
\end{rema}

\begin{figure}[htb]

\psfrag{p1}{$p_1$}
\psfrag{p2}{$p_2$}
\psfrag{s1}{$s_1$}
\psfrag{q1}{$q_1$}
\psfrag{q2}{$q_2$}
\psfrag{l1u}{$\ell_{1}^u$}
\psfrag{l2u}{$\ell_{2}^u$}
\psfrag{y1u}{$y_{1}^u$}
\psfrag{y2u}{$y_{2}^u$}
\psfrag{G1}{$\GG^{u}$}
\psfrag{Ss}{$\SS^{s}$}

   \includegraphics[width=9cm]{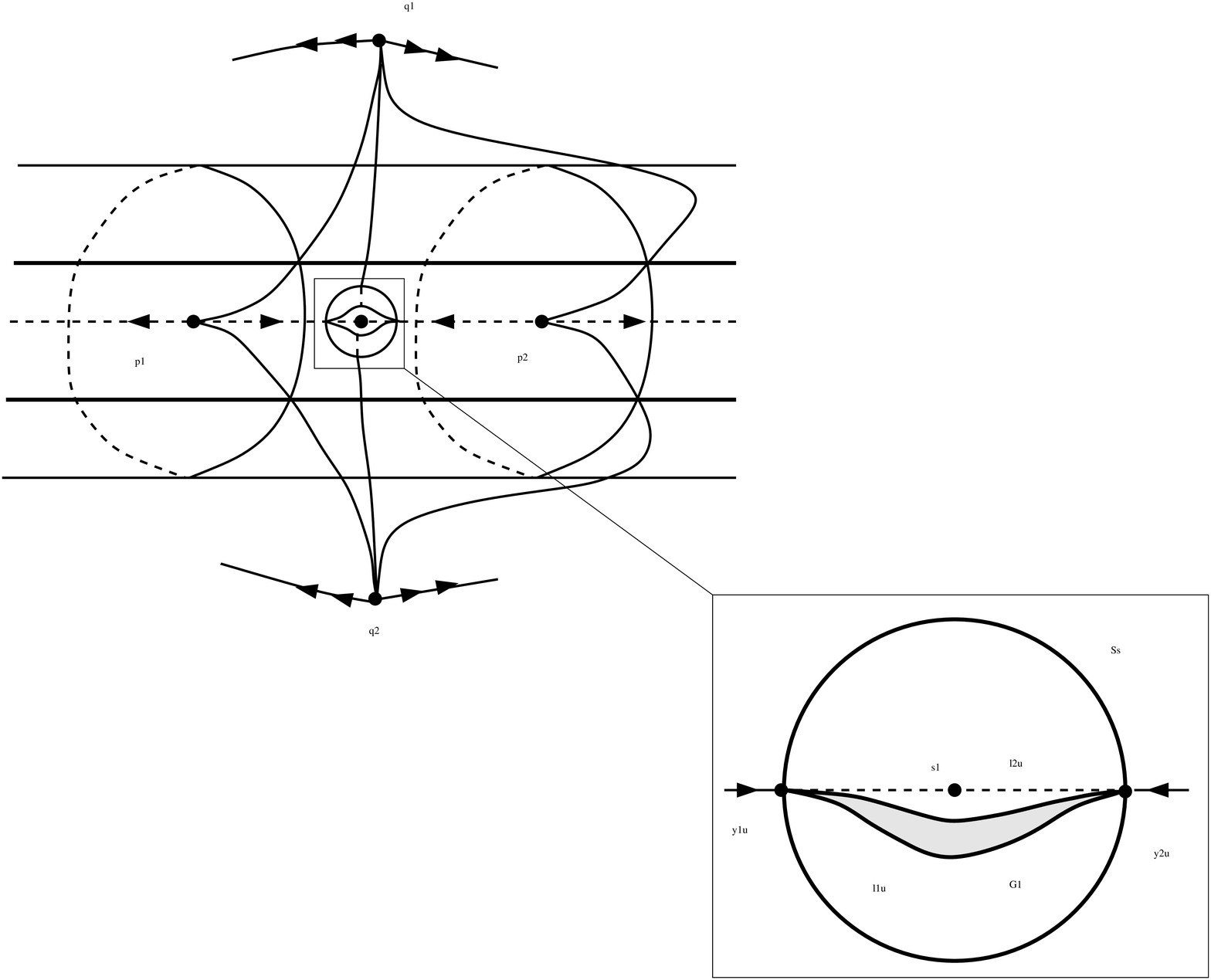}\hspace{1cm}
    \caption{The set $\GG^u$.}
    \label{f.Gu}
\end{figure}

Consider the sets $\GG^u$ and $\GG^s$ (the set $\GG^ u$ is depicted
in Figure~\ref{f.Gu}),
$$
\GG^u \eqdef \overline{W^u(r_1) \cap \SS^s} \quad \mbox{and} \quad
\GG^s\eqdef \overline{W^s(s_1) \cap \SS^u}.
$$

\begin{lemm}\label{l.banana}
The set  $\GG^u$ is a topological two-disk bounded by
$ \ell_1^u \cup \ell_2^u \cup \{y_1^u\} \cup \{y_2^u\}$, where
$$
\ell_1^u \subset W^u(q_1) \quad  \mbox{and} \quad \ell_2^u \subset
W^u(q_2)
$$
are disjoint (open) simple curves whose endpoints are $y_1^u$ and
$y_2^u$. The closed curves $\overline{\ell_1^u}$ and
$\overline{\ell_2^u}$ have the same tangent direction  at their
endpoints $y^u_1$ and $y^u_2$. The disk $\GG^u$ has two cusps at
the points
$y^u_1$ and $y^u_2$.

Similarly,
the set $\GG^s$ is a topological two-disk bounded by $\ell_1^s \cup
\ell_2^s \cup \{y_1^s\} \cup \{y_2^s\}$, where
$$
\ell_1^s \subset W^s(p_1) \quad \mbox{and} \quad \ell_2^s \subset
W^s(p_2)
$$
are disjoint (open) simple curves whose endpoints are $y_1^s$ and
$y_2^s$. The closed curves $\overline{\ell_1^s}$ and
$\overline{\ell_2^s}$ have the same tangent direction  at their
endpoints $y^s_1$ and $y^s_2$.
 The disk $\GG^s$ has two cusps at $y_1^s$ and $y^s_2$.
\end{lemm}

We consider the following notation, given an interval $[t_1,t_2]$
and a set $A$ we let
$$
X_{[t_1,t_2]}(A)\eqdef \bigcup_{t\in [t_1,t_2]} X_t(A).
$$

\begin{proof} We only prove the lemma for the set $\GG^u$, the proof
for $\GG^s$ is identical.

Note that the sphere $\SS^s$ intersects every orbit of the set
$W^s(s_1)\setminus \{s_1\}$ in exactly one point.
Thus, since
 the rectangle $R$ in Remark~\ref{r.stablecylinders} is the closure
of $W^u(r_1)\cap W^s(s_1)\cap \TT^2$, the
positive orbit of any point in the interior of $R$ intersects
$\SS^s$ in exactly one point. Hence the set $\GG^u$ is the
closure of the ``projection"  along the orbits of $X$ of the
interior of $R$ into  $\SS^s$, that is,
$$
\GG^u= \overline{
 X_{[0,\infty)} \big( \mbox{int}(R) \big) \cap \SS^s}.
$$
By construction the set $\GG^u$ is a topological two-disk. We next
describe its boundary.

By equation \eqref{e.boundaryR} the  boundary of $R$ consists of
the segments $a_i^s\subset W^s(p_i)$ and $b_i^u\subset W^u(q_j)$,
$i=1,2$. Furthermore, the interiors of the segments $b_1^u$ and
$b_2^u$  are contained in $W^s(s_1)$. Denote by $\ell^u_1$ and
$\ell^u_2$ the ``projections'' by the flow of $X$ of these interiors
into $\SS^s$, that is,
$$
\ell_i^u \eqdef X_{[0,\infty)} \big( \mbox{int} (b_i^u) \big) \cap
\SS^s.
$$

Consider any sequence $(x_n)$ of points in the interior of $R$
accumulating to the side $a_i^s$ of $R$. Note that  the (positive)
orbit of $x_n$ by the flow of $X$ goes arbitrarily close to the
saddle singularity $p_i$ before intersecting
$\SS^s$ at a point $y_n$. By construction, the sequence $(y_n)$ converges
to $y_i^u=W^u(p_i)\cap \SS^s$. Indeed, for any given curve $b\subset R$ transverse to  $X$ joining the
sides $a_1^s$ and $a_2^s$ of $R$ the intersection of the sphere $\SS^s$ 
and the positive orbit of $b$ by the flow 
 $X$ (i.e., the
``projection" of  $b$ into $\SS^s$ by the flow)
 is a
curve $\ell_b$ joining $y^u_2$ and $y_2^u$ (these points are in
the closure of $\ell_b$). In particular, $y_1^u$ and $y_1^u$ are
the endpoints of $\ell^u_i$, $i=1,2$. 

Bearing in mind equation~\eqref{e.sides}
and the definitions of $y_{i}^u$, $\ell^u_{i}$, and $\Ga(q_i)$,
$i=1,2$, we get the following:
\[
\begin{split}
 \overline{\ell^u_i} &= \ell^u_i
\cup \{y_1^u,y_2^u\}= W^u(\Ga(q_i)) \cap \SS^s,
\\
\partial \GG^u&= \big( W^u(\Ga(q_1)) \cap \SS^s \big) \cup \big( W^u(\Ga(q_2)) \cap \SS^s
\big)= \ell^u_1 \cup \ell^u_2  \cup \{y_1^u,y_2^u\}.
\end{split}
\]
 This completes the description of the set $\partial
\GG^u$.

It remains to see that $y_1^u$ and $y_2^u$ are cusps of $\GG^u$.
By Lemma~\ref{l.tangent}
the curves $\overline{\ell^u_1}$ and
$\overline{\ell^u_2}$ are tangent at $y_i^u$ to 
$E^{cs}(y_i^u)\cap T_{y_i^u}(\SS^s)$. 
To see that the point $y_i^u$ is a cusp of $\GG^u$
 it is enough to note that the interior of $\GG^u$ is disjoint from the circle
 $\eta^s$.
Thus the disk $\GG^u$ is the ``thin component" of $\SS^s\setminus
\partial \GG^u$. This ends the proof of the lemma. 
\end{proof}

\begin{rema}
\label{r.lastinclusion}
With the notations above, the following inclusions hold
$$
(\SS^s\setminus \GG^u) \subset W^u(r_2) \quad \mbox{and} \quad
(\SS^u\setminus \GG^s\subset W^s(s_2)).
$$
\end{rema}






\subsection{The diffeomorphism time-one map $X_1$}\label{ss.gluing}
Let $X_1$ denote the time-one map  of the vector field $X$ and
define the diffeomorphism $F_0\eqdef X_1.$ Note that $F_0$ is a
Morse-Smale diffeomorphism whose non-wandering set consists of the
sinks $s_1$ and $s_2$, the saddles of $s$-index two $p_1$ and
$p_2$, the saddles of $s$-index one $q_1$ and $q_2$, and the
sources $r_1$ and $r_2$. Note that the invariant manifolds of these
points for the vector field $X$ and for the diffeomorphism $F_0$ coincide. We only write
$W^i(x,X)$ or $W^i(x,F_0)$ to emphasize the role of $X$ or $F_0$,
otherwise we just write $W^i(x)$.

\begin{figure}[htb]

\psfrag{l1u}{$\ell_{1}^u$}
\psfrag{l2u}{$\ell_{2}^u$}
\psfrag{y1u}{$y_{1}^u$}
\psfrag{y2u}{$y_{2}^u$}
\psfrag{L1u}{$L_{1}^u$}
\psfrag{L2u}{$L_{2}^u$}
\psfrag{Y1u}{$Y_{1}^u$}
\psfrag{Y2u}{$Y_{2}^u$}
\psfrag{Gu}{$\GG^{u}$}
\psfrag{FGu}{$F_0(\GG^{u})$}
\psfrag{Eu}{$\EE^{u}$}
\psfrag{G1}{$\GG^{u}$}
\psfrag{Ss}{$\SS^{s}$}
\psfrag{Ds}{$\Delta^{s}$}

   \includegraphics[width=9cm]{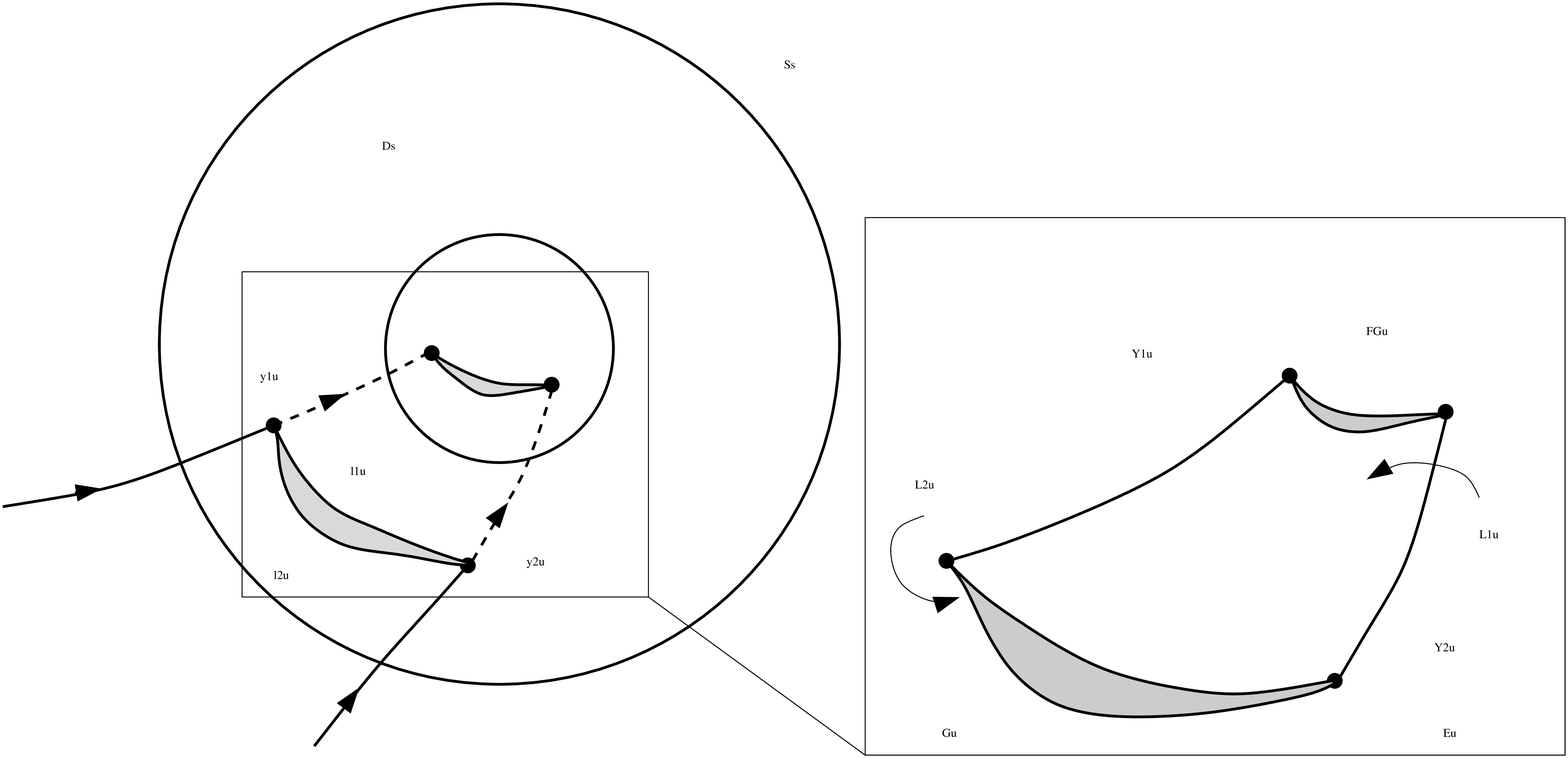}\hspace{1cm}
    \caption{The sets $\EE^u, \LL_i^u,$ and $\YY^u_i$.}
    \label{f.Eu}
\end{figure}

Consider the fundamental domain $\De^s$ of $W^s(s_1)$ for $F_0$
bounded by $\SS^s$ and $F_0(\SS^s)$. Note that
$$
\De^s \eqdef X_{[0,1]} (\SS^s) = \BB^s\setminus
\mbox{int}(F_0(\BB^s)) \simeq \SS^s\times [0,1].
$$
Let
$$
\EE^u\eqdef X_{[0,1]} (\GG^u), \quad \LL^u_i \eqdef X_{[0,1]}
(\ell^u_i), \quad \YY^u_i \eqdef X_{[0,1]} (y^u_i).
$$
These sets are depicted in
Figure~\ref{f.Eu}.

By construction, the set $\EE^u$ is a topological ball bounded by the disks $\GG^u$ and
$F_0(\GG^u)$, the ``rectangles" $\LL_1^u$ and $\LL^u_2$, and
the curves $\YY^u_1$ and $\YY^u_2$. By the definition of $\GG^u$
and of the cylinders $\CC_1^u$ and $\CC^u_2$ and by
Lemma~\ref{l.banana} (see also Remark~\ref{r.lastinclusion}), we have that
\begin{equation} \label{e.inEs} 
\overline{W^u (r_1) \cap
\De^s}=\EE^u \quad \mbox{and} \quad \De^s\setminus \EE^u \subset
 W^u(r_2).
\end{equation}
Since $\ell^u_i\subset W^u(q_i,X)$ and $y^u_i\subset W^u(p_i,X)$
we have
\begin{equation}\label{e.LsYs}
\LL^u_i \subset W^u(q_i,F_0) \quad \mbox{and} \quad
 \YY^u_i \subset W^u(p_i,F_0).
\end{equation}
Note that the common boundary of the rectangles $\LL^u_1$ and
$\LL^u_2$ are the curves $\YY_1^u$ and $\YY^u_2$. Moreover,
the closure of $\LL_i^u$ is  tangent to the center unstable bundle
$E^{cu}$ along the curves $\YY_1^u$ and $\YY^u_2$.

We similarly define the fundamental domain $\De^u$ of $W^u(r_1)$
for $F_0$ bounded by $\SS^u$ and $F_0(\SS^u)$. As before
$$
\De^u \eqdef X_{[0,1]} (\SS^u) \simeq \SS^u\times [0,1].
$$
We let
$$
\EE^s\eqdef X_{[0,1]} (\GG^s), \quad \LL^s_i \eqdef X_{[0,1]}
(\ell^s_i), \quad \YY^s_i \eqdef X_{[0,1]} (y^s_i)
$$
and as above we have
\begin{equation}\label{e.LuYu}
\LL^s_i \subset W^s(p_i) \quad \mbox{and} \quad
 \YY^s_i \subset W^s(q_1).
\end{equation}
Moreover, the ``rectangles" $\LL^s_1$ and $\LL^s_2$  are tangent
at the curves  $\YY_1^s$ and $\YY^s_2$ to the center stable bundle
$E^{cs}$. By construction, $\EE^s$ is bounded by the disks $\GG^s$
and $F_0(\GG^s)$, the ``rectangles" $\LL_1^s$ and $\LL^s_s$, and
the curves $\YY^s_1$ and $\YY^s_2$. Finally,
\begin{equation}
\label{e.inEu} \overline{W^s (s_1) \cap \De^u}=\EE^s \quad
\mbox{and} \quad \De^u\setminus \EE^s \subset
 W^s(s_2).
\end{equation}

\section{A diffeomorphism on $\SS^2\times \SS^1$ obtained by a
surgery}\label{s.surgery}

\subsection{The surgery}\label{ss.surgery}
 In this section 
identify some regions of $\SS^3$ by a local
diffeomorphism $\Psi$. This surgery  provides
 a diffeomorphism $F_\Psi$ defined on $\SS^2 \times \SS^1$ 
 induced by $F_0=X_1$ and the quotient of $\SS^3$ by $\Psi$.
We will see in Section~\ref{s.psi} 
that for an appropriate choice of $\Psi$ the diffeomorphism $F_\Psi$ 
has fragile cycles.

With the notation in Section~\ref{s.auxiliary},
 consider a diffeomorphism $\Psi\colon \De^s\to \De^u$
such that
\begin{itemize}
\item
$\Psi(\SS^s)= \SS^u$ and  $\Psi(F_0(\SS^s))=F_0(\SS^u)$,
\item
$\Psi\circ F_0|_{\SS^s}= F_0|_{\SS^u}\circ \Psi$,
\item
$D\Psi\circ DF_0|_{\SS^s}= DF_0|_{\SS^u}\circ D\Psi$.
\end{itemize}

Recall that $\BB^s\subset W^s(s_1)$ and $\BB^u\subset W^u(r_1)$
are small balls containing $s_1$ and $r_1$, respectively.
In the set $ \SS^3\setminus \Big( \mbox{int}
\big( F_0(\BB^s) \big)\cup \mbox{int}\big( \BB^u \big) \Big)$ we 
 identify the
points $x\in \De^s$ and $\Psi(x)\in \De^u$ obtaining the
quotient space
$$
M\eqdef \left( \SS^3\setminus \Big( \mbox{int} \big( F_0(\BB^s)
\big)\cup \mbox{int}\big( \BB^u) \big)\Big )\right)  \Big\slash
 \Psi.
$$
The set $M$ is a $C^1$-manifold diffeomorphic to $\SS^2\times \SS^1$
 and the diffeomorphism $F_0$ induces  a diffeomorphism $F_\Psi
\colon M\to M$. 

Denote by $\pi$ the projection $\pi\colon  \SS^3\setminus \Big(
\mbox{int} \big( F_0(\BB^s) \big)\cup \mbox{int}\big( \BB^u)
\big)\Big ) \to M$ that associates to $x$ its class
$\pi(x)=[x]$ by the equivalence relation induced by $\Psi$. For
notational simplicity, if $x\not\in
\De^s \cup \De^u$ we simply write $x$ instead of $\pi(x)$.
Write
$$
\De\eqdef \big(\De^s\cup \De^u \big)/\Psi=\pi(\De^s)=\pi(\De^u).
$$
In what follows we write $\KK^u_i(F_\Psi)\eqdef \pi (\KK^u_i)$,
where $\KK=\YY,\LL,\GG,\EE$ and $i=,1,2, \emptyset$.

\begin{rema}\label{r.delta} {\em{The set
$$
\De\setminus \big (\LL_1^u (F_\Psi) \cup \LL_2^u  (F_\Psi)\cup
\YY^u_1 (F_\Psi) \cup \YY^u_2 (F_\Psi) \cup \GG^u (F_\Psi) \cup
 F_\Psi ( \GG^u (F_\Psi) ) \big)
$$
has two connected components. The set $\EE^u (F_\Psi)$ is the
component that is diffeomomorphic to $\GG^u \times [0,1]$. This set is a topological
ball.
 There is a similar characterization for
$\EE^s(F_\Psi)$. }}
\end{rema}

\subsection{The dynamics of $F_\Psi$}
Note that the orbits of $F_\Psi$ disjoint from $\De$ are the
projection of the orbits of the diffeomorphism $F_0$.
Thus $s_2$ is a sink and $r_2$ is a source of $F_\Psi$
(note that the sink $s_1$ and the source $r_1$ are removed
in the construction of $M$). Similarly, the points $p_1$ and $p_2$ are saddles
of $s$-index $2$ of $F_\Psi$ and the points $q_1$ and $q_2$ are saddles of $s$-index
$1$ of $F_\Psi$. Observe that $F_{\Psi}$ 
can  have further periodic points, but by 
 Remark~\ref{r.time}
these points
have period larger
than $10$.

Using the identification by $\Psi$ and the properties of $F_0$,
 we
have the following characterization of the sets $\LL^{u,s}_i(F_\Psi)$
and  $\YY^{u,s}_i(F_\Psi)$, $i=1,2$:
\begin{equation}\label{e.inclusions}
\begin{split}
\LL^u_i(F_\Psi)&=\{x\in \De \, \colon \, x \in W^u(q_i,F_\Psi)
\,\, \mbox{and}\,\, F_\Psi^{-j}(x) \notin \De \,\, \mbox{for
all}\,\, j\ge
2\},\\
\YY^u_i(F_\Psi)&=\{x\in \De \, \colon \, x \in W^u(p_i,F_\Psi)
\,\, \mbox{and}\,\, F_\Psi^{-j}(x) \notin \De \,\, \mbox{for
all}\,\, j\ge 2\},
\\
\LL^s_i(F_\Psi)&=\{x\in \De \, \colon \, x \in W^s(p_i,F_\Psi)
\,\, \mbox{and}\,\, F_\Psi^{j}(x) \notin \De \,\, \mbox{for
all}\,\, j\ge 2\},
\\
\YY^s_i(F_\Psi)&=\{x\in \De \, \colon \, x \in W^s(q_i,F_\Psi)
\,\, \mbox{and}\,\, F_\Psi^{j}(x) \notin \De \,\, \mbox{for
all}\,\, j\ge 2\}.
\end{split}
\end{equation}

\begin{rema}\label{r.paredesFPsi}
{\em{
The normally hyperbolic curves $\Ga_{i,j}$ of $F_0$ (recall Remark~\ref{r.Gaij}) 
are disjoint from $\BB^s\cup \BB^u$, thus their projections on $M$ (also denoted by $\Ga_{i,j}$)
are normally hyperbolic curves of $F_\Psi$ .
Observe also that by construction,
the interior of $\Ga_{i,j}$ is contained in 
$W^s( p_i,F_\Psi) \pitchfork W^u(q_j,F_\Psi)$, recall equation
\eqref{e.curvasGa}.

We continue to use the notation 
$\Ga(q_j)=\Ga_{1,j}\cup\Ga_{2,j}$ and $\Ga(p_i)=\Ga_{i,1}\cup \Ga_{i,2}$.
}}
\end{rema}

With the previous notation we have  that
$$
\LL^u_j(F_\Psi) \cup \YY^u_1 (F_\Psi) \cup \YY^u_2 (F_\Psi)= \overline{\LL^u_j(F_\Psi)}
$$
is a connected component of $W^u(\Ga (q_j), F_\Psi)\cap \De$.

\begin{lemm}[Invariant manifolds and their intersections]
\label{l.bananasfirst} Consider $x\in \De$.
\begin{enumerate}
\item \label{banana1}
If $x\notin \EE^u (F_\Psi)$ then  $x\in W^u(r_2,F_\Psi)$ and thus
it is not chain recurrent,
\item \label{banana2}
if $x\notin  \EE^s (F_\Psi)$ then $x\in W^s(s_2,F_\Psi)$ and thus
it is not chain recurrent.
\item \label{banana3}
if $x\in \LL^u_i (F_\Psi)$ then $x\in W^u(q_i,F_\Psi)$,
\item
\label{banana4} if
 $x\in \YY^u_i (F_\Psi)$ then $x\in W^u(p_i,F_\Psi)$,
\item
\label{banana5} if $x\in \LL^s_i (F_\Psi)$ then $x\in
W^s(p_i,F_\Psi)$, and
\item
\label{banana6} if
 $x\in  \YY^s_i (F_\Psi)$ then $x\in
W^s(q_i,F_\Psi)$,
\end{enumerate}
\end{lemm}

\begin{proof}The first  item follows immediately from equation \eqref{e.inEs} and the definition of $F_\Psi$.
Similarly, the second  item follows  from \eqref{e.inEu}.
Items (\ref{banana3})-(\ref{banana6}) follow from equation \eqref{e.inclusions}.
\end{proof}

\begin{figure}[htb]

\psfrag{z}{$z$}
\psfrag{w}{$w$}
\psfrag{Y1s}{$\YY_{1}^s(F_\Psi)$}
\psfrag{L1s}{$\LL_{1}^s(F_\Psi)$}
\psfrag{Y1u}{$\YY_{1}^u (F_\Psi)$}
\psfrag{L1u}{$\LL_{1}^u(F_\Psi)$}
\psfrag{Eu}{$\EE^{u}(F_\Psi)$}
\psfrag{Es}{$\EE^{s}(F_\Psi)$}
\psfrag{Ss}{$\SS^{s}$}
\psfrag{homoclinic}{$\De^u\cap R(F_\Psi)$}

   \includegraphics[width=9cm]{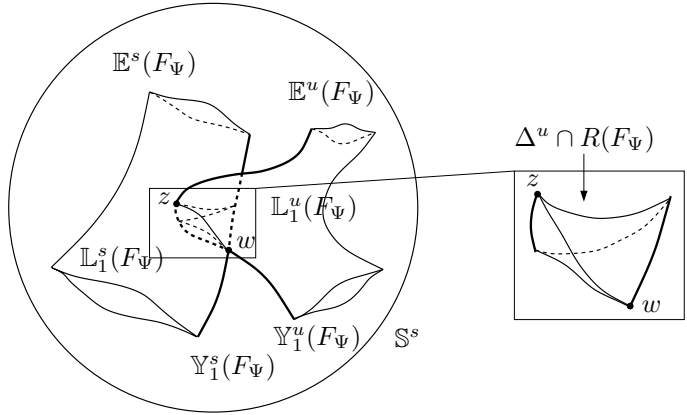}\hspace{1cm}
    \caption{The recurrence region of $F_\Psi$.}
    \label{f.lemma}
\end{figure}

An immediate consequence of the first two items of Lemma~\ref{l.bananasfirst}
is the following:

\begin{coro}\label{c.inextremis}
Let $x\in \De$ be a chain recurrent point for $F_\Psi$. Then $x\in
\EE^s (F_\Psi)\cap  \EE^u (F_\Psi)$.
\end{coro}

A straightforward  consequence of Lemma~\ref{l.bananasfirst}
and Remark~\ref{r.paredesFPsi}
 is the following:

\begin{coro}\label{c.bananas} Consider $x\in \De$.
\begin{enumerate}
\item \label{cbanana1}
If $x\in \LL_i^u (F_\Psi)\cap \LL_j^s (F_\Psi)$ then $x\in
W^u(q_i,F_\Psi) \pitchfork W^s(p_j,F_\Psi)$.
\item \label{cbanana2}
If $x\in \LL_i^u (F_\Psi) \cap \YY_j^s (F_\Psi)$ then $x\in
W^u(q_i,F_\Psi) \cap W^s(q_j,F_\Psi)$.
\item \label{cbanana3}
If $x\in \YY_i^u (F_\Psi) \cap \LL_j^s (F_\Psi)$ then $x\in
W^u(p_i,F_\Psi) \cap W^s(p_j,F_\Psi)$.
\item \label{cbanana4}
If $x\in \YY_i^u (F_\Psi) \cap \YY_j^s (F_\Psi)$ then $x\in
W^u(p_i,F_\Psi) \cap W^s(q_j,F_\Psi)$.
\item
\label{cbanana5} The interior of $\Ga_{i,j}$ is contained in
$W^s(p_i,F_\Psi)\pitchfork W^u(q_j,F_\Psi)$.
\end{enumerate}
\end{coro}


Recall that a periodic point is called {\emph{isolated}} if its chain recurrent class coincides
with its (finite) orbit.

\begin{lemm}[Heterodimensional cycles and trivial homoclinic classes]\label{l.hetero}
$\,$
\begin{enumerate}
\item \label{hh1}
If  $\YY^u_1 (F_\Psi)\cap \YY^s_1 (F_\Psi)\ne\emptyset$  then
$F_\Psi$ has a heterodimensional cycle associated to $p_1$ and
$q_1$.
\item \label{hh2}
If $\YY^u_1 (F_\Psi)\pitchfork \LL^s_1 (F_\Psi)\neq \emptyset$
then the homoclinic class of $q_1$ is non trivial.
\item \label{hh3}
If $\YY^s_1 (F_\Psi)\pitchfork  \LL^u_1 (F_\Psi)\neq \emptyset$
then homoclinic class of $p_1$ is non trivial.
\item
\label{hh4} If  $\YY^u_1 (F_\Psi)\cap \EE^s (F_\Psi) =
\emptyset$ then  $p_1$ is an isolated saddle. 
\item
\label{hh5}
 If  $\YY^s_1 (F_\Psi) \cap \EE^u (F_\Psi) = \emptyset$ then
$q_1$ is  an isolated saddle.
\end{enumerate}
\end{lemm}
\begin{proof}
To prove the first item just recall that by item (\ref{cbanana5})
in Corollary~\ref{c.bananas} the interior of $\Ga_{1,1}$ is
contained in $W^s(p_1,F_\Psi)\pitchfork W^u(q_1,F_\Psi)$. By item
(\ref{cbanana4}) in Corollary~\ref{c.bananas}, if
$\YY^u_1(F_\Psi)\cap \YY^s_1 (F_\Psi)\ne \emptyset$ then $W^u(p_1,F_\Psi)\cap
W^s(q_1,F_\Psi)\ne\emptyset$ and thus there is a heterodimensional cycle
associated to $p_1$ and $q_1$.

Items \eqref{hh2} and \eqref{hh3} follow from Lemma~\ref{l.bananasfirst}.

To prove item \eqref{hh4} 
we use  the following simple fact whose proof we omit.

\begin{rema}\label{r.nonisolated}{\em{
Let $p$ be  a hyperbolic saddle that is non-isolated. 
Then its stable/unstable manifold contains points of its chain recurrence class 
that do not belong to its orbit. In particular, if a hyperbolic fixed point  $p$ is such that
$W^u(p)\setminus \{p\}$
(respectively, $W^s(p)\setminus \{p\}$)
 is contained in the stable (resp. unstable) manifolds of some sinks
 (resp. sources) then it  is isolated.}}
\end{rema}

Note that every point $x\in W^u
(p_1,F_\Psi)$, $x\ne p_1$, in the separatrix of $W^u(p_1,F_\Psi)$
that does not contain $y_1^u$ is contained in $W^s(s_2,F_\Psi)$.
Thus this separatrix does not contain chain recurrent points.
Hence it is enough to consider points in the  separatrix of
$W^u(p_1,F_{\Psi})$ containing $y_1^u$.  Note that 
$W^u(p_1,F_{\Psi})\cap \De$ contains a fundamental domain of 
$W^u(p_1,F_{\Psi})$ that is contained
in $\YY^u_1(F_\Psi)$.
Thus, by Remark~\ref{r.nonisolated}, if $p_1$ is not isolated then 
$\YY^u_1(F_\Psi)$  must contain some point of the chain recurrence class of $p_1$.
Thus in such a case $\YY^u_1(F_\Psi)$ cannot be contained in $W^s(s_2,F_\Psi)$.

Suppose that  $\YY^u_1(F_\Psi)\cap
\EE^s (F_\Psi)= \emptyset$. Then, by item \eqref{banana2} in
Lemma~\ref{l.bananasfirst}, one has that
$\YY^u_1(F_\Psi)\subset W^s(s_2)$. By the above discussion
this implies that $p_1$ is isolated, 
proving item
\eqref{hh4}.

The proof of item  \eqref{hh5} is identical to the previous one and thus it is omitted.
\end{proof}

\section{Dynamics in a neighborhood of $F_\Psi$}
\label{ss.neighborhood} 
In this section we consider
diffeomorphisms $F$ in a small $C^1$-neighborhood of $F_\Psi$.
The hyperbolic-like properties of the objects
introduced in Section~\ref{s.surgery} for the diffeomorphism $F_\Psi$
allows us to define their continuations for nearby diffeomorphisms
and thus to repeat these constructions.
In particular, the continuations of the hyperbolic  points
$s_2,r_2,p_1,p_2,q_1$, and $q_2$  of $F_\Psi$ are defined. We will omit the
dependence on $F$ of these continuations.

As the arguments in this section are similar to those in 
Sections~\ref{s.auxiliary} and \ref{s.surgery},
some of these constructions will be just sketched. We now go to the details of our
constructions.

Consider the spheres $\pi(\SS^s)$ and $F(\pi(\SS^s))$ and
denote by $\De_F$ the closure of the connected component of
$M\setminus \big( \pi(\SS^s) \cup F(\pi(\SS^s)) \big)$ which  is
close to $\De$ (i.e., the component that is in the same side of
$\pi(\SS^s)$ as $\De$). The set $\De_F$ is diffeomorphic to
$\SS^2\times [0,1]$ and varies continuously with $F$ in the
$C^1$-topology. In particular,  by Remark~\ref{r.time}, if $x\in \De_F$ and
$F^i(x)\in \De_F$ then $|i|\ge 9$.

Bearing in mind the definitions of
the sets $\LL^{s,u}_{i}(F_\Psi)$ and $\YY^{s,u}(F_\Psi)$
in \eqref{e.inclusions},
we define their ``continuations"  $\LL^{s,u}_{i}(F)$ and $\YY^{s,u}(F)$
for $F$ close to $F_\Psi$ by
\[
\begin{split}
\LL^u_i(F)&=\{x\in \De_F \, \colon \, x \in W^u(q_i,F) \,\,
\mbox{and}\,\, F^{-i}(x) \notin \De_F \,\, \mbox{for all}\,\, i\ge
2\},\\
\YY^u_i(F)&=\{x\in \De_F \, \colon \, x \in W^u(p_i,F) \,\,
\mbox{and}\,\, F^{-i}(x) \notin \De_F \,\, \mbox{for all}\,\, i\ge
2\},
\\
\LL^s_i(F)&=\{x\in \De_F \, \colon \, x \in W^s(p_i,F) \,\,
\mbox{and}\,\, F_\Psi^{i}(x) \notin \De_F \,\, \mbox{for all}\,\,
i\ge 2\},
\\
\YY^s_i(F)&=\{x\in \De_F \, \colon \, x \in W^s(q_i,F) \,\,
\mbox{and}\,\, F^{i}(x) \notin \De_F \,\, \mbox{for all}\,\, i\ge
2\}.
\end{split}
\]

\begin{rema}\label{r.continuousinextremis}
{\em{The sets $\LL^{s,u}_i(F)$ and $\YY^{s,u}_i(F)$, $i=1,2$,
depend continuously on $F$.}}
\end{rema}

Note  that  the closed curves  $\Ga_{i,j}$ are normally hyperbolic 
for $F_\Psi$ (recall Remark~\ref {r.paredesFPsi}). 
Thus for every $F$ close to $F_\Psi$ there
are defined their continuations, denoted by
$\Ga_{i,j}(F)$, that depend continuously on $F$. 
These curves join the saddles $p_i$ and $q_j$ and
their interiors are center stable manifolds of $p_i$ and center
unstable manifolds of $q_j$. Finally, from the normal
hyperbolicity of $\Ga_{i,j}(F)$, compact parts of the invariant
manifolds $W^s(\Ga_{i,j}(F))$ and $W^u(\Ga_{i,j}(F))$ depend
continuously on $F$.

Observe that $W^u(q_i,F)\setminus W^{uu}(q_i,F)$ (resp.
$W^s(p_i,F)\setminus W^{ss}(p_i,F)$) has  two connected
components (separatrices), denoted by $W^{u,+}(q_i,F)$ and
$W^{u,-}(q_i,F)$ (resp. $W^{s,\pm}(p_i,F)$). We
choose these components such that the following holds:

\begin{rema}[Invariant manifolds of $\Ga_{i,j}(F)$]
\label{r.compatibleintersections}
$\,${\em{
\begin{itemize}
\item 
$W^u(\Ga_{1,j}(F))\setminus W^u(p_1,F)=  W^u(\Ga_{2,j}(F))\setminus W^u(p_2,F)= W^{u,+}(q_j,F)\cup W^{uu}(q_j,F)$.
\item 
$W^s(\Ga_{i,1}(F))\setminus W^s(q_1,F) =  W^s(\Ga_{i,2}(F))\setminus W^s(q_2,F)= W^{s,+}(p_i,F)\cup$\newline
$W^{ss}(p_i,F)$.
\item $W^u(\Ga_{i,1}(F))\cap W^u(\Ga_{i,2}(F))= W^u(p_i,F)$.
\item $W^s(\Ga_{1,j}(F))\cap W^s(\Ga_{2,j}(F))= W^s(q_j,F)$.
\end{itemize}}}
\end{rema}

By Lemma~\ref{l.centerbundles} and using the notation in Remark~\ref{r.notation}, for every $F$ close
to $F_\Psi$ there is a unique invariant central bundle $E^{cs}_F$ 
(resp. $E^{cu}_F$)
defined on
$W^u(p_i,F)$ (resp. $W^s(q_j,F)$) and transverse to the strong stable (resp. unstable) direction at $p_i$
(resp. $q_j$).
The central bundles $E^{cs}_F$ and $E^{cu}_F$ 
depend continuously on $F$.

\begin{rema}\label{r.compatiblecentral}
{{\em 
The manifolds with boundary $W^u(\Ga_{i,1}(F))$ and
$W^u(\Ga_{i,2}(F))$ are tangent along $W^u(p_i,F)$ to the plane
field $E^{cs}_F$. As in Remark~\ref{r.semiplanes},  for $x\in
W^u(p_i,F)$, the vectors of $E^{cs}_F$ entering in the interior of
$W^u(\Ga_{i,j}(F))$ define an half plane $E^{cs}_{+,F}(x)$.

Analogously,
the manifolds with boundary $W^s(\Ga_{1,j}(F))$ and
$W^u(\Ga_{2,j}(F))$ are tangent along $W^s(q_j,F)$ to the plane
field $E^{cu}_F$. For $x\in W^s(p_j,F)$, the vectors of $E^{cu}_F$
entering in the interior of $W^(\Ga_{i,j}(F))$ define a half plane
$E^{cu}_{+,F}$.
}}
\end{rema}

As in \eqref{e.Gaqi} and \eqref{e.Gapi}, we
define the sets
$$
\Ga(q_j,F)\eqdef \Ga_{1,j}(F)\cup\Ga_{2,j}(F) \quad \mbox{and}
\quad \Ga(p_i,F)\eqdef \Ga_{i,1}(F)\cup\Ga_{i,2}(F).
$$
Then
$$
W^u(\Ga(q_j,F))= W^u(\Ga_{1,j}(F))\cup W^u(\Ga_{2,j}(F))
$$
is a $C^1$-surface with boundary
whose compact parts
 depend continuously on $F$. Moreover,
$$
W^u(\Ga(q_1,F))\cap W^u(\Ga(q_2,F))= W^u(p_1,F)\cup W^u(p_2,F)
$$
 The surfaces  $W^u(\Ga(q_1,F))$ and $W^u(\Ga(q_2,F))$
 depend continuously on $F$ and
  are
tangent to $E^{cs}_F$ along this intersection. This last assertion is just a
version of Lemma~\ref{l.tangent} for $F$ close to $F_\Psi$.

Similarly, the
compact parts of the surfaces $W^s(\Ga(p_1,F))$ and $W^s(\Ga(p_2,F))$ depend
continuously on $F$ and they are tangent to $E^{cu}_F$
along their intersection $W^s(q_1,F)\cup W^s(q_2,F)$.

Using  the previous notation  we get that
\[
\begin{split}
&\overline{\LL^u_j(F)}=\LL^u_j(F)\cup \YY^u_1(F)\cup \YY^u_2(F),\\
& \overline{\LL^s_i(F)}=\LL^s_i(F)\cup \YY^s_1(F)\cup \YY^s_2(F).
\end{split}
\]
The set $\overline{\LL^u_j(F)}$ is the connected
component of $W^u(\Ga(q_j,F))\cap\De_F$ whose negative iterates
$F^{-i}(\overline{\LL^u_j(F)})$, $i\ge 2$,
 are disjoint from $\De_F$. Similarly, the set  $\overline{\LL^s_i(F)}$ is the connected
component of $W^u(\Ga(p_i,F)\cap\De_F$ whose positive iterates
larger than $2$ are disjoint from $\De_F$.

As a consequence of the previous constructions, we get

\begin{lemm}\label{l.rectanglescontinuous} 
The sets
$\overline{\LL^u_j(F)}$ and $\overline{\LL^s_i(F)}$ are ``rectangles" 
depending  continuously on $F$ (for the $C^1$ topology).
\end{lemm}

The sets $\GG^{s,u} (F)$ and $\EE^{s,u} (F)$ are defined similarly as in the case
$F_\Psi$. The set $\GG^u(F)$ is the topological disk with two cusps
(these cuspidal points are  in $W^u(p_1,F)$ and $W^u(p_2,F)$) whose boundary
is the union of $W^u(\Ga(q_1,F))\cap \De_F$ and   
$W^u(\Ga(q_2,F))\cap \De_F$. There is an analogous definition for $\GG^s(F)$.
Note that by construction these sets depend continuously on $F$.

Finally, the set $\EE^u(F)$ is the topological ball bounded by $\overline{\LL^u_1(F)}$, $\overline{\LL^u_2(F)}$, $\GG^u(F)$ and $F(\GG^u(F))$ that is close
to $\EE^u(F_\Psi)$. There is a similar definition for the set $\EE^s(F)$.
By construction the sets $\EE^{s,u}(F)$ depend continuously on $F$.

There is the following reformulation of
Corollary~\ref{c.inextremis} and Lemma~\ref{l.hetero} for 
diffeomorphisms $F$ close to $F_\Psi$.

\begin{lemm} 
\label{l.inextremis}
Consider a diffeomorphism $F$ close to $F_\Psi$.
\begin{enumerate}
\item \label{F1}
If $x\in \De_F$ is a chain recurrent point for $F$ then $x\in
\EE^s(F) \cap \EE^u(F)$.
\item \label{F2}
If  $\YY^u_1(F)\cap \YY^s_1(F)\ne\emptyset$  then $F$ has a
heterodimensional cycle associated to $p_1$ and $q_1$.
\item \label{F3}
If $\YY^u_1(F) \pitchfork \LL^s_1 (F)\neq \emptyset$ then the
homoclinic class of $q_1$ is non trivial.
\item \label{F4}
If $\YY^s_1(F) \pitchfork  \LL^u_1 (F) \neq \emptyset$ then
homoclinic class of $p_1$ is non trivial.\item
\label{F5} If  $\YY^u_1 (F) \cap \EE^s (F) = \emptyset$ then
 $p_1$ is  isolated.
\item
\label{F6}
 If  $\YY^s_1 (F) \cap \EE^u (F) = \emptyset$ then 
$q_1$ is  isolated.
\end{enumerate}
\end{lemm}

Observe also that by construction
 $$
 \De_F\subset \EE^s(F)
 \subset
 W^s(s_2,F)
\quad \mbox{and} \quad 
 \De_F\setminus \EE^u(F)
 \subset
 W^u(r_2,F).
 $$
As in Corollary~\ref{c.inextremis}, 
a consequence of these inclusions is the following.

\begin{lemm}
For every $F$ close to $F_\psi$, every point of $\De_F$  that is chain recurrent is contained in $\EE^s(F)\cap \EE^u(F)$.
\end{lemm}

\section{Choice of the local diffeomorphism $\Psi$} \label{s.psi} 

Lemmas~\ref{l.bananasfirst} means that,
 for the diffeomorphism  $F_\Psi$, 
the existence 
of heterodimensional cycles and homoclinic
intersections for $p_1,p_2,q_1$, and $q_2$ 
 depend on the intersections of the sets $\LL^{u,s}_i
(F_\Psi)$, $\YY^{u,s}_i (F_\Psi)$, $\EE^u (F_\Psi)$, and
$\EE^s(F_\Psi)$. The choice of the identification map
$\Psi$ determines these intersections. Lemma~\ref{l.inextremis}
explains how these properties are translated for diffeomorphisms
close to $F_\Psi$.

\begin{figure}[htb]

\psfrag{r}{$r$}
\psfrag{Y1s}{$\YY_{1}^s (F_\Psi)$}
\psfrag{Y1u}{$\YY_{1}^u (F_\Psi)$}
\psfrag{Eu}{$\EE^{u}(F_\Psi)$}
\psfrag{Es}{$\EE^{s}(F_\Psi)$}
\psfrag{Ecu+}{$E^{cu}_+$}
\psfrag{Ecs+}{$E^{cs}_+$}
   \includegraphics[width=5cm]{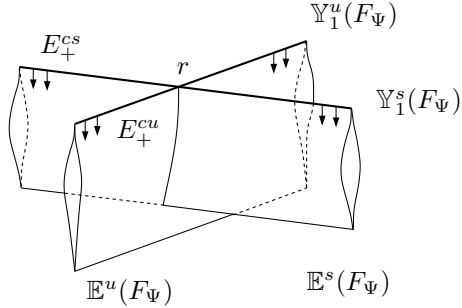}\hspace{1cm}
    \caption{Conditions T and D}
    \label{f.td}
\end{figure}

We assume that the local diffeomorphism $\Psi$ is such that 
the diffeomorphism $F_\Psi$
satisfies the following two conditions, see Figure~\ref{f.td}:

\smallskip

\noindent {\bf{(T) Topological hypothesis:}} There is a point
$z\in \De$ such that
$$\
\YY^s_1(F_\Psi)\cap \YY^u_1(F_\Psi)=
\YY^s_1(F_\Psi)\cap \EE^u(F_\Psi)=\YY^u_1(F_\Psi)\cap  
\EE^s(F_\Psi)= \{z\}.
$$

\noindent {\bf{ (D) Differentiable hypothesis:}}
\begin{itemize}
\item
The intersection of the semi-planes $E^{cs}_{+}(z)\cap
E^{cu}_{+}(z)$ is a half straight line. This implies that the
intersection $\YY^u_1(F_\Psi)\cap \YY^s_1(F_\Psi)$ at $z$ is quasi-transverse.
\item
In a neighborhood of $z$, the sets $\EE^s(F_\Psi)$ and
$\EE^u(F_\Psi)$ are locally in the same side of any locally
defined surface containing the curves $\YY^u_1(F_\Psi)$ and
$\YY^s_1(F_\Psi)$.
\end{itemize}

\begin{prop}\label{p.subsigma}
Suppose that that $F_\Psi$
satisfies conditions (T) and (D). Then there
are a $C^1$-neighborhood $\cU_{F_\Psi}$ of $F_\Psi$ and a
codimension one submanifold $\Si$ containing $F_\Psi$
such that:
\begin{enumerate}
\item
The set $\cU_{F_\Psi}\setminus \Si$ is the union of two
disjoint open sets $\cU^+_{\Si}$ and $\cU^-_{\Si}$  such
that:
\begin{itemize}
\item
for every $G\in \cU^+_{\Si}$ the saddle $p_{1}$ is isolated
and the homoclinic class of $q_{1}$ is non-trivial, and
\item
for every $G\in \cU^-_{\Si}$ the saddle $q_1$ is isolated
and the homoclinic class of $p_1$ is non-trivial.
\end{itemize}
\item
Every diffeomorphism $G\in \Si$ has a heterodimensional
cycle associated to $p_{1}$ and $q_{1}$.
\end{enumerate}
\end{prop}

We have the following corollary:

\begin{coro}\label{co.fragile}
The submanifold $\Si$ consists of diffeomorphisms $F$
having fragile cycles associated to $p_1$ and $q_1$.
\end{coro}

Note that Theorem~\ref{t.fragil} follows immediately from Proposition~\ref{p.subsigma} and Corollary~\ref{co.fragile}.

Corollary~\ref{co.fragile} is a consequence of Proposition~\ref{p.subsigma}
and the following simple fact about cycles and chain recurrence  
classes:

\begin{rema}[Cycles and chain recurrence classes]\label{r.cycleschains}
{\em{Let $F$ be a diffeomorphism with a heterodimensional cycle
associated to two transitive hyperbolic sets $L$ and $K$.
Then every pair of saddles $p\in L$ and $q\in K$ are in the same chain
recurrence class of $F$. In particular, both saddles $p$ and $q$ are not isolated.}}
\end{rema}

\begin{proof}[Proof of Corollary~\ref{co.fragile}]
We argue by contradiction. If there is $F\in \Si$ such
that the cycle is not fragile then there is a diffeomorphism $G$ close to $F$ with a
robust cycle associated to a pair of transitive hyperbolic sets
$L\ni p_1$ and $K\ni q_1$. Then by 
Remark~\ref{r.cycleschains} the saddles $p_1$ and $q_1$ 
are not isolated for $G$. Finally, as the cycle is robust, we can assume that $G\not\in
\Si$, that is, $G\in \cU_\Si^+\cup \cU_\Si^-$. Since $p_1$ is isolated if $G\in \cU^+_{\Si}$
this implies that $G\not\in \cU^+_{\Si}$. Similarly, as $q_1$
is isolated if $G\in \cU^-_{\Si}$ we have that $G\not\in
\cU^-_{\Si}$. This contradicts the fact that the cycle associated to $G$ is robust.
\end{proof}

\subsection{Proof of Proposition~\ref{p.subsigma}}
To define the submanifold $\Sigma$
we take the unitary vector $\overrightarrow n$ normal to the plane
$T_z(\YY^u_{1} (F_\Psi))\oplus T_z(\YY^s_{2}(F_\Psi))$  pointing
to the ``opposite" direction of $\EE^u(F_\Psi)$ and
$\EE^s(F_\Psi)$. Fix small $\epsilon>0$ and for $F$ close to $F_\Psi$
consider the family of one-dimensional disks
$$
\left\{
\YY^s_{1}(F)+t\, \overrightarrow{n}\right\}_{t\in [-\epsilon,
\epsilon]}.
$$

The main step of the proof of the proposition if the following lemma whose
proof we postpone.

\begin{lemm}\label{l.preprop}
There is a small
neighborhood $\cU_{F_\Psi}$ of $F_{\psi}$ such that for every $F\in \cU_{F_\Psi}$
there is a unique parameter $t=\tau_F$, depending continuously on $F$,
such that 
$$
\big( \YY^s_1(F)+\tau_F\,
\overrightarrow{n}\big) \cap \YY^u_2(F)\ne\emptyset.
$$ 
There are the following three possibilities according to the value of
$\tau_F$:
\begin{itemize}
\item If $\tau_F=0$ then the diffeomorphism $F$ has a heterodimensional
cycle associated to $p_1$ and $q_1$.
\item If $\tau_F>0$ then $p_{1}$  is isolated
and the homoclinic class of $q_1$ is non-trivial.
\item If $\tau_F<0$ then  $q_{1}$ is isolated 
and the homoclinic class of $p_1$ is
non-trivial.
\end{itemize}
\end{lemm}

\begin{figure}[htb]

\psfrag{z}{$z$}
\psfrag{Y1s}{$\YY_{1}^s(F)$}
\psfrag{Y1u}{$\YY_{1}^u (F)$}
\psfrag{text1}{$q_1$ is isolated}
\psfrag{text2}{$H(p_1)$ is non-trivial}
\psfrag{text3}{$p_1$ is isolated}
\psfrag{text4}{$H(q_1)$ is non-trivial}

 \includegraphics[width=8cm]{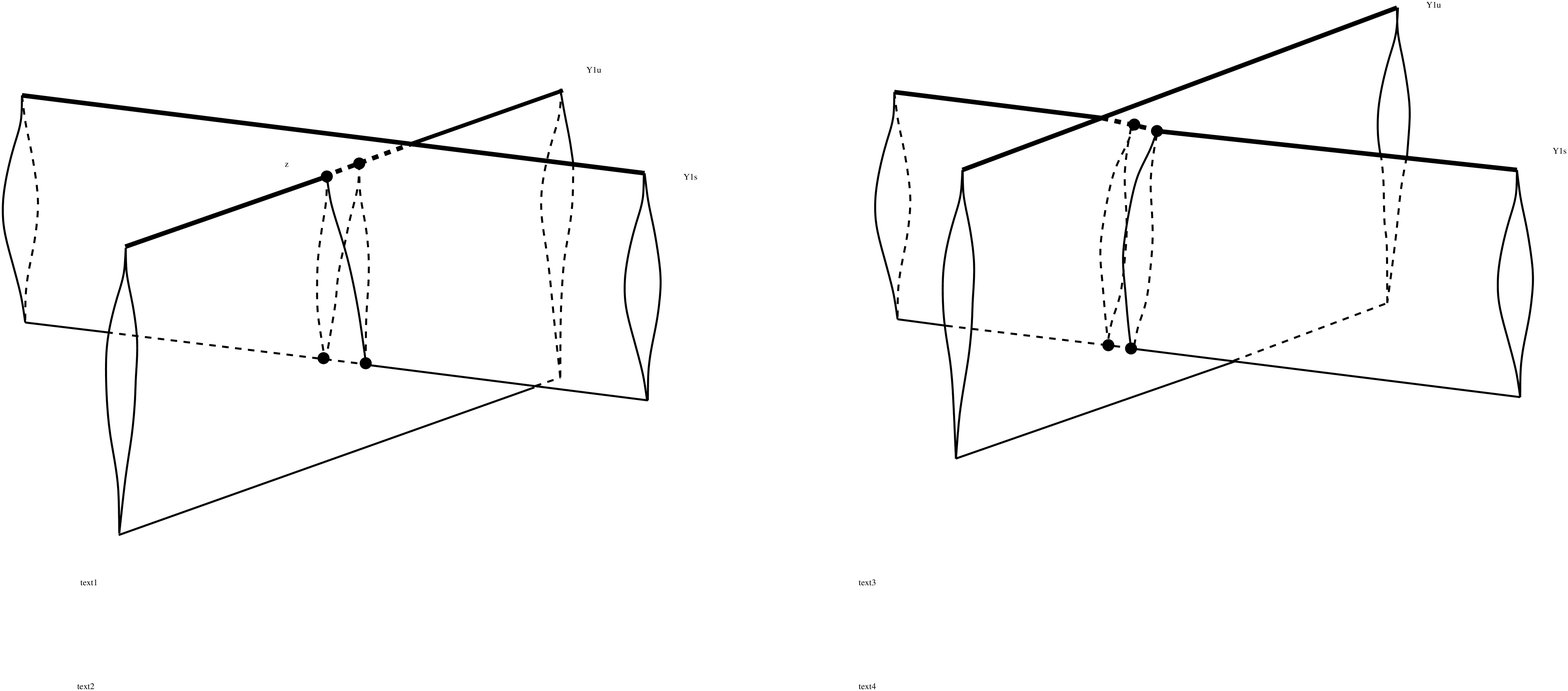}\hspace{1cm}
    \caption{Lemma~\ref{l.preprop}}
    \label{f.isolated}
\end{figure}

\begin{proof}[Proof of Proposition~\ref{p.subsigma}]
In view of Lemma~\ref{l.preprop} we let
$$
\Sigma\eqdef \{ F\in \cU_{F_\Psi}\, \mbox{such that $\tau_F=0$}\}.
$$
Then the set $\cU_{F_\Psi} \setminus \Si$
 has two components  $\cU_{\Si}^+$ and $\cU_\Si^-$.
 The component $\cU_{\Si}^+$
  consists of the
 diffeomorphisms $F$ such that $\tau_F>0$ and
 $\cU_{\Si}^-$ consists of the
 diffeomorphisms $F$ with $\tau_F<0$. The proposition now follows immediately from
 Lemma~\ref{l.preprop}.
 \end{proof}

\begin{proof}
Consider the surface
$$
\cY^s_1(F)=\bigcup_{t\in [-\ve,\ve]} \YY^s_{1}(F)+t\, \overrightarrow{n}.
$$
Recall that $\YY^s_1(F)$ and $\YY^u_1(F)$ depend
continuously on $F$ (see Remark~\ref{r.continuousinextremis}). 
Thus the surface $\cY^s_1(F)$ also depends
continuously on $F$. Note that, by hypothesis, $\cY^s_1(F_\Psi)$
is transverse to $\YY^u_1(F_\Psi)$ and this intersection is 
 the point $z$ in condition (T). Thus for $F$ close to $F_\Psi$ the
intersection $\cY^s_1(F) \cap\YY^u_{1} (F)$ also consists of
exactly one point $z_F$ depending continuously on $F$. Thus there is 
exactly one parameter $\tau_F$ with
$$
z_F \in \YY^s_1(F) + \tau_F\, \overrightarrow n.
$$
Moreover, the parameter 
$\tau_F$
depends continuously on $F$.

If $\tau_F=0$ then $\YY^u_1(F)\cap \YY^u_2(F)
\ne\emptyset$. Thus by item (\ref{F2}) in Lemma~\ref{l.inextremis}
the diffeomorphism $F$ has a heterodimensional cycles associated
to $p_1$ and $q_1$.

The choice of $\overrightarrow n$ implies that  for $\tau_F>0$ the sets
$\YY^u_{1}(F)$ and $\EE^s(F)$ are disjoint. Thus
by item (\ref{F5}) in Lemma~\ref{l.inextremis}
the point
 $p_{1}$ is isolated for $F$. 
 
Similarly, for $\tau_F>0$ one has that $\YY^s_1(F)$ intersects
$\EE^u(F)$. Thus, if $\tau_F$ is small enough, one has that 
$\YY^s_1(F)\pitchfork
\LL^u_1(F)$. Hence by item (\ref{F3}) in Lemma~\ref{l.inextremis}
the homoclinic class $H(q_1,F)$ is
non-trivial.

A similar argument shows that for $\tau_F<0$ the sets $\YY^u_{1}(F)$
and $\EE^u(F)$ are disjoint. As above item (\ref{F6}) 
in Lemma~\ref{l.inextremis}
implies that 
the point $q_1$ is isolated. Also for 
$\tau_F<0$ one has  that $\YY^u_1(F)\pitchfork \LL^s_1(F)$
and item (\ref{F4}) in Lemma~\ref{l.inextremis}
implies that the homoclinic class $H(p_1,F)$ is non-trivial. 
This completes the proof of the lemma.
\end{proof}

\begin{rema}{\em{
Our construction can be done such that the saddles $q_2$ and $q_1$ are homoclinically related for every small  $t>0$. For that it is enough to choose the
diffeomorphisms $\Psi$ such that $\YY^s_2(F_\Psi)\pitchfork \LL^u_1(F_\Psi)$,
see Corollary~\ref{c.bananas}.

Observe that in this case one has $H(q_2,F_t)=H(q_1,F_t)$ for $t>0$. 
However, for $t=0$ the saddle $q_1$ escapes from the non-trivial homoclinic
class of $q_2$. Surprisingly, in such a case one can generate robust cycles associated to $q_2$ and $p_1$ (this follows from \cite{BDK}) but not associated to $q_1$ and $p_1$.
}}
\label{r.recall}
\end{rema}

\section{Discussion} \label{s.conclusion}

Our construction provides examples of fragile cycles relating two
saddles  $p_1$ and $q_1$  of different indices (for simplicity we will
omit the dependence on the diffeomorphisms). This
construction has a ``prescribed" part concerning the
relative positions of the invariant manifolds of these saddles. But
this prescribed dynamics involves only the ``cuspidal" regions
$\YY^s_1$ and $\YY^u_1$ of 
$\EE^s$ and $\EE^u$, respectively. The rest of the intersection $\EE^s\cap \EE^u$ can
be chosen arbitrarily. This provides a
lot of ``freedom"
for the global dynamics, for instance, for the
behavior of the other
saddles of the diffeomorphism and their invariant manifolds
(recall Remark~\ref{r.recall}). 
Hence, without further assumptions on the dynamics outside the
cycle the global dynamics cannot be described.  

\subsection{Partial hyperbolicity, wild dynamics, and
fragile cycles}
 A first ingredient of our construction is the
gluing map $\Psi$ that plays a key role for determining the
resulting global dynamics. 

\begin{ques}\label{q.ph}
Can the resulting dynamics be partially hyperbolic? More
precisely, does it exist a gluing map $\Psi$ so that
the chain recurrence class of $p_1$ for $F_\Psi$ has a partially hyperbolic
splitting with three $1$-dimensional bundles?
\end{ques}

We expect a positive answer to this question. This would
imply that the phenomena associated to these fragile cycles 
would also occur in the most rigid non-hyperbolic setting of partial
hyperbolicity with one-dimensional center.

\smallskip

We next discuss how fragile cycles can be involved in 
the generation of {\emph{wild dynamics}}
(roughly, persistent coexistence of infinitely
many homoclinic classes, see \cite{bible} for further details). 

First note that perturbations of the diffeomorphism $F$ in
 the fragile cycles submanifold
 $\Sigma$ in Theorem~\ref{t.fragil} expels  periodic
points in the chain recurrence class that 
simultaneously
contains $p_1$ and $q_1$:
the saddle $p_1$ is expelled by diffeomorphisms in the component $\cU^+_{\Si}$ and the saddle  $q_1$ by diffeomorphisms in the component $\cU^-_{\Si}$. 
Thus the chain recurrence class fall into pieces.
On the other hand, 
there is a natural question about whether fragile cycles may generate a ``cascade" of fragile cycles:

\begin{ques}\label{q.repeat}
Can the submanifold of fragile cycles $\Si$ be accumulated by (co-dimension one) submanifolds  consisting of fragile cycles?
\end{ques}
 
If the answer to this question would be positive then a ``fragile cycles
configuration"  could  be repeated  generating infinitely many different
chain recurrence classes. 
Being very optimistic, positive answers to the two questions above
could provide the first examples  of
wild dynamics in the partially hyperbolic (with
one-dimensional center) setting. 
We will discuss further questions about wild dynamics later. 

\subsection{Chain recurrence classes}\label{ss.chain}
Consider a diffeomorphism $F$ in the submanifold $\Si$ of fragile
cycles.  By construction, the intersection
$W^s(p_1)\pitchfork W^u(q_1)$ contains a curve $\Ga_{1,1}$ joining
$p_1$ to $q_1$. Moreover, the intersection $W^u(p_1)\cap W^u(q_1)$
is exactly the orbit of a quasi-transverse heteroclinic point
$x_F$. This cyclic configuration implies that the chain recurrence
classes of $p_1$ and $q_1$  coincide and contain the curve $\Ga_{1,1}$
and the orbit of $x_F$, 
recall also Remark~\ref{r.cycleschains}.

\begin{figure}[htb]

\psfrag{x}{$x$}
\psfrag{y}{$y$}
\psfrag{r}{$r$}
\psfrag{p}{$p_1$}
\psfrag{q}{$q_1$}
\psfrag{g11}{$\Ga_{1,1}$}

   \includegraphics[width=4.2cm]{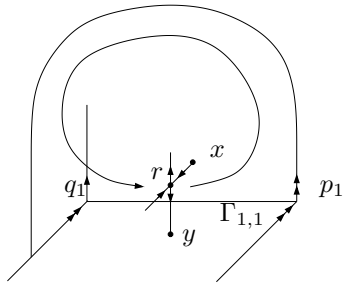}\hspace{1cm}
    \caption{Periodic points in $C(p_1)=C(q_1)$.}
    \label{f.twisted}
\end{figure}

 We can  now perform 
perturbations of $F\in \Si$ preserving the cycle, that is, the
resulting diffeomorphisms continue to belong to $\Si$. In this way, and using
for instance 
the arguments in \cite{BD08}, one can slightly modify the central
eigenvalues of the saddles $p_1$ and $q_1$ in the cycle 
(corresponding to the tangent
direction of the connection $\Ga_{1,1}$) 
to generate ``new" periodic saddles $r$ whose
orbits pass arbitrarily close to $p_1$ and $q_1$ and belong
to the chain recurrence class of $C(p_1)=C(q_1)$. The latter fact is a consequence of the geometry of the cycle that guarantees that indeed $W^u(r)\cap W^s(p_1)\ne \emptyset$ and $W^s(r)\cap W^u(q_1)\ne\emptyset$, see Figure~\ref{f.twisted}.

\begin{ques}\label{q.Si} For $C^1$-generic diffeomorphisms
$F\in\Si$, is the heteroclinic curve $\Ga_{1,1}$ contained in the
closure of a set of periodic points?
Can the periodic points in that set  be chosen in the same
homoclinic class or (even more) homoclinically related?
\end{ques}

Concerning this  question,  observe that our construction
provides saddles that are in the same recurrence class (the one of
$p_1$ and $q_1$). Note also that by \cite{BC04}, for
$C^1$-generic diffeomorphisms, the chain recurrence classes of
periodic points always coincide with their homoclinic classes. The
difficulty here is that we consider $C^1$-generic diffeomorphisms
in the codimension one submanifold $\Si$ (which is a meager set). 
Thus the result in \cite{BC04}
cannot be applied. Indeed, for $F\in \Si$ the homoclinic classes of
$p_1$ and $q_1$ are both trivial and
hence different from their non-trivial chain
recurrence classes since $p_1,q_1\in C(p_1)=C(q_1)$.

We believe that the answer to Question~\ref{q.Si} is positive.
In such a case, it would be  interesting to consider the central Lyapunov exponents of the periodic orbits in the chain recurrence class $C(p_1)=C(q_1)$.
This is our next topic.

\subsection{Lyapunov exponents} The papers
\cite{DHRS09,LOR,DG} consider examples of diffeomorphisms having twisted
heterodimensional cycles (recall Section~\ref{s.introduction}) and analyze the spectrum of central Lyapunov exponents of the (non-trivial) homoclinic classes involved in the cycle. In the examples in \cite{LOR,DG} this spectrum has a {\emph{gap}}. 
The existence of this gap is related to the fact that the homoclinic
class $H(p)$ considered contains a saddle $q$ of different index
that looks like a ``cuspidal corner point"  (exactly as the points $p_1$ and
$q_1$ of our construction)
and satisfies $q\in H(p)$ but $H(q)=\{q\}$.
So the saddle $q$ is topologically an extreme point but also dynamically is also an extreme point as it is not homoclinically related to other saddles in the class $H(p)$.
We believe that precisely this is reflected by the Lyapunov
spectrum that has  one gap.


The previous comments and the fact that 
the chain recurrence classes considered in this paper have ``two cuspidal corner
points" (the saddles $p_1$ and $q_1$) that are also dynamically extremal lead to the following question.

\begin{ques} Does there exist a
 gluing map $\Psi$ such that there is a neighborhood $\cU$ of
$F_\Psi$ such that ``persistently" in $\cU\cap \Si$ the
diffeomorphisms have  two (or more) gaps in the spectrum of central Lyapunov
exponents of the periodic orbits in the chain recurrence class $C(p_1)$? 
Here the term ``persistent" is purposely vague and may mean $C^1$-generic or $C^1$-dense, for instance.
\end{ques}

Observe that by \cite{ABCDW07} for a $C^1$-generic diffeomorphism the spectrum of Lyapunov exponents 
of homoclinic classes has no gaps.
Hence, one would need to consider diffeomorphisms of non-generic type.

\subsection{Collision, collapse, and birth of classes}
Consider an arc of diffeomorphisms $(F_t)_{t\in[-1,1]}$ that intersects $\Si$  transversely at $t=0$. Assume that $F_t\in \cU^+_\Si$ for $t>0$ and 
$F_t\in \cU^-_\Si$  for $t<0$. 
Our construction implies that for $t>0$ the class $H(q_1,F_t)$ is non-trivial, 
$H(p_1,F_t)=\{p_1\}$, and both classes are disjoint. Similarly, for $t<0$ the class $H(p_1,F_t)$ is non-trivial, $H(q_1,F_t)=\{q_1\}$, and both classes are disjoint. Since $F_0\in \Si$ and $H(q_1,F_0)=\{q_1\}$ and 
$H(p_1,F_0)=\{p_1\}$, each of these classes collapses to a single point at $t=0$. Similarly (or symmetrically), for $t=0$ the chain recurrence classes of $p_1$ and $q_1$ collide at $t=0$ and contain the heteroclinic segment $\Ga_{1,1}$, when one of these classes collapses to a point for $t\neq 0$. This illustrates the lower semi-continuous dependence of homoclinic classes and the upper semi-continuity of chain recurrence classes.
 
 Let us observe that \cite{BCGP} provides a  locally      $C^1$-dense set of diffeomorphisms where homoclinic classes are properly contained
 in a robustly isolated chain recurrence class. This construction is somewhat 
  similar to the one in this paper 
and involves  a heterodimensional cycle relating cuspidal corner points of the class.

Coming back to our construction, it would be interesting to understand how for $t>0$ the points of  $H(q_1,F_t)$ escape from the class or ``disappear" as $t\to0^+$. As we discussed above, in some cases the saddle $q_1$ is accumulated by saddles of the same index when $t=0$. These saddles are not homoclinically related to $q_1$ but they are in its chain recurrence class. This indicates that there are (infinitely many) saddles that ``escape" from the homoclinic class of $q_1$ but not from the chain recurrence class of $q_1$ as $t$ evolves. Though we do not know how these saddles are homoclinically related. 
Note that there is completely symmetric scenery for the saddle $p_1$.

\begin{ques}
Is there a diffeomorphism $F\in\Si$ with  infinitely many different
homoclinic classes, all of them contained in the chain recurrence class of $p_1$ and $q_1$?
\end{ques}

In view of the previous discussion, a simpler question is if the diffeomorphisms in $\Si$  can be chosen having wild  dynamics close to the cycle. More precisely, 

\begin{ques}
Does there exist $F\in\Si$ with infinitely many different chain recurrence classes accumulating to $p_1$ or to $q_1$?
\end{ques}

Finally, let us observe that the previous setting is somehow reminiscent to the setting of the geometrical Lorenz attractor: at the point of ``bifurcation", when a singular cycle occurs, infinitely many
orbits of the vector field transform into heteroclinic orbits of a singularity, see for
instance \cite{BLMP}.


\section*{Acknowledgments}
This paper was partially supported by  CNPq, FAPERJ, and Pronex (Brazil),
Agreement France-Brazil in Mathematics, and ANR Project DynNonHyp 
BLAN08-2$_-$313375.
LJD thanks the warm hospitality of Institut de Math\`ematiques de Bourgogne.




\begin{thebibliography}{100}


\bibitem{A03}
F. Abdenur, \emph{Generic robustness of spectral decompositions,}
Ann. Sci. \'Ecole Norm. Sup. {\bf 36} (2003), 213--224.

\bibitem{ABCDW07} F.Abdenur, Ch. Bonatti, S. Crovisier, L. J.
D\'{\i}az, and L. Wen, {\emph{Periodic points and homoclinic classes,}}
Ergod. Th. \& Dynam. Sys. \textbf{27} (2007), 1--22



\bibitem{BLMP}
R. Bam\'on, R. Labarca, R. Ma\~n\'e, and M. J. Pacifico,
{\emph{The explosion of singular cycles,}}
Publ. Math. I.H.\'E.S. {\bf 78} (1993), 208--231. 


\bibitem{ASY06}
K. T. Alligood, E. Sander, and J.A. Yorke, {\emph{
Three-dimensional crisis: Crossing bifurcations and unstable
dimension variability,}} Phys. Rev. Lett. \textbf{96} (2006),
244103.




\bibitem{bible}
Ch. Bonatti, \emph{Towards a global view of dynamical systems, for
the $C^1$-topology,} preprint Institut de Math\'ematiques
de Bourgogne (2010) and to appear in Ergod. Th. \& Dynam.
Sys..


\bibitem{BC04}
Ch. Bonatti and S. Crovisier, {\emph{R\'ecurrence et
g\'en\'ericit\'e,}} Inventiones Math.  {\bf 158} (2004), 33--104.

\bibitem{BCGP}
Ch. Bonatti, S. Crovisier, N. Gourmelon, and R. Portrie,
{\emph{Tame non-robustly transitive diffeomorphisms,}} in preparation.





\bibitem{BD08}
Ch.~Bonatti and L. J. D\'iaz, {\emph{Robust heterodimensional
cycles and $C^{1}$-generic dynamics,}} Jour. Inst. of
Math. Jussieu \textbf{7}  (2008),  469--525




\bibitem{BDK}
Ch.~Bonatti, L. J. D\'iaz, and S. Kiriki, {\emph{Stabilization of heterodimensional cycles,}} in preparation.

\bibitem{BDVbook}
 Ch.~Bonatti, L. J. D\'iaz, and M. Viana,
 {\emph{Dynamics beyond uniform hyperbolicity,}}
 Encyclopaedia of Mathematical Sciences (Mathematical Physics) \textbf{102},
 Mathematical physics, III,
\textit{Springer Verlag}, 2005.


\bibitem{Con}
C. Conley, {\emph{Isolated invariant sets and Morse index,}} 
CBMS Regional Conf. Ser. Math. {\bf 38}, Providence (RI), Am. Math. Soc. 1978. 











\bibitem{DG}
    L.~J.~D\'iaz and K. Gelfert, {\emph{Porcupine-like horseshoes: Transitivity, Lyapunov spectrum, and phase transitions,}}
 {\tt{arXiv:1011.6294}}.

\bibitem{DHRS09}
    L.~J.~D\'iaz, V.~Horita, I.~Rios, and M.~Sambarino,
    \emph{Destroying horseshoes via heterodimensional cycles: generating bifurcations inside homoclinic classes},
Ergod. Th. \& Dynam.
Sys. ~\textbf{29} (2009), 433--473.











\bibitem{LOR}
    R.~Leplaideur, K.~Oliveira, and I.~Rios,
    \emph{Equilibrium states for partially hyperbolic horseshoes},
   Ergod. Th. \& Dynam.
Sys., to appear.
   
   

\bibitem{Gugu}
C. Moreira, {\emph{There are no $C^1$-stable intersections of regular Cantor sets,}}
preprint IMPA, {\tt{http://w3.impa.br/$\sim$gugu/}}.

 \bibitem{Newhouse}
S. Newhouse, 
{\emph{Hyperbolic limit sets,}}
Trans. Amer. Math. Soc. \textbf{167} (1972), 125--150.

\bibitem{N68}
S. Newhouse, {\emph{Nondensity of axiom ${\rm A}({\rm a})$ on
$S\sp{2}$,}} Global Analysis (Proc. Sympos. Pure Math., Vol. XIV,
Berkeley, Calif., Amer. Math. Soc. (1968),  191--202.



 \bibitem{N79}
S.\,E. Newhouse, 
{\emph{The abundance of wild hyperbolic sets and
non-smooth stable sets for diffeomorphisms,}} 
Publ.\ Math.\
I.H.\'E.S.
 \textbf{50} (1979), 101--151






















\end{thebibliography}
\end{document}